\documentclass[11pt,letter]{article}
\usepackage{amsfonts}
\usepackage{amssymb}
\usepackage{amsmath}
\usepackage{amsthm}
\usepackage{subfigure}
\usepackage{psfrag}
\usepackage{cite}
\usepackage{url}

\usepackage[dvips]{graphicx}

\title{A Geometrical Study of Matching Pursuit Parametrization}
\date{\today}

\usepackage{color}

\newcommand{\cqfd}{$\hfill\blacksquare$}
\renewcommand{\vec}[1]{{\bf #1}}

\newcommand{\ud}{\mathrm{d}} 
\newcommand{\Nbb}{\mathbb{N}}
\newcommand{\Zbb}{\mathbb{Z}}

\newcommand{\Rbb}{\mathbb{R}}
\newcommand{\Cbb}{\mathbb{C}}
\newcommand{\scp}[2]{\langle #1, #2 \rangle}
\newcommand{\bscp}[2]{\big\langle #1, #2 \big\rangle}
\newcommand{\inv}[1]{\frac{1}{#1}}
\newcommand{\tinv}[1]{{\textstyle\frac{1}{#1}}}
\newcommand{\dict}{{\rm dict}}

\newtheorem{theorem}{Theorem}

\newtheorem{lemma}{Lemma}
\newtheorem{proposition}{Proposition}
\newtheorem{corollary}{Corollary}
\newtheorem{definition}{Definition}
\newtheorem{question}{Question}

\begin{document}

\maketitle
\begin{abstract}
This paper studies the effect of discretizing the parametrization of a dictionary used for Matching Pursuit decompositions of signals. Our approach relies on viewing the continuously parametrized dictionary as an embedded manifold in the signal space on which the tools of differential (Riemannian) geometry can be applied. The main contribution of this paper is twofold. First, we prove that if a discrete dictionary reaches a minimal density criterion, then the corresponding discrete MP (dMP) is equivalent in terms of convergence to a weakened hypothetical continuous MP. Interestingly, the corresponding weakness factor depends on a density measure of the discrete dictionary. Second, we show that the insertion of a simple geometric gradient ascent optimization on the atom dMP selection maintains the previous comparison but with a weakness factor at least two times closer to unity than without optimization. Finally, we present numerical experiments confirming our theoretical predictions for decomposition of signals and images on regular discretizations of dictionary parametrizations.
\end{abstract}

\section{Introduction}
\label{sec:intro}

There has been a large effort in the last decade to develop analysis techniques that decompose non-stationary signals into elementary components, called {\it atoms}, that characterize their salient features~\cite{basispursuit,maza93,greedisgood,bandelets,uncertainty}. In particular, the matching pursuit (MP) algorithm has been extensively studied~\cite{maza93,NeffVideoMP,vetterli,stochasticMP,harmonicaudio,expoconvergence,aprioriMP} to expand a signal over a redundant dictionary of elementary atoms, based on a greedy process that selects the elementary function that best matches the residual signal at each iteration. Hence, MP progressively isolates the structures of the signal that are coherent with respect to the chosen dictionary, and provides an adaptive signal representation in which the more significant coefficients are first extracted. The progressive nature of MP is a key issue for adaptive and scalable communication applications~\cite{frossard3D,inloopquantization}.

A majority of works that have considered MP for practical signal approximation and compression define the dictionary based on the discretization of a parametrized prototype function, typically a scaled/mo\-du\-la\-ted Gaussian function or its second derivative~\cite{NeffVideoMP,gribonval2001fmp,Rosa}. 
An orthogonal 1-D or 2-D wavelet basis is also a trivial example of such a discretization even if in that case MP is not required to find signal coefficients; a simple wavelet decomposition is computationally more efficient. Works that do not directly rely on a prototype function either approximate such a parametrized dictionary based on computationally efficient cascades of filters~\cite{SubbandDicoCDV,FastMPBull,NeffDico}, or attempt to adapt a set of parametrized dictionary elements to a set of training signal samples based on vector quantization techniques~\cite{gainshape,SchmidSaugeon}. Thus, most earlier works define their dictionary by discretizing, directly or indirectly, the parameters of a prototype function.

The key question is then: \emph{how should the continuous parameter space be discretized ?} A fine discretization results in a large dictionary which approximates signals efficiently with few atoms, but costs both in terms of computational complexity and atom index entropy coding. Previous works have studied this trade-off empirically~\cite{Rosa,NeffVideoMP}. In contrast, our paper focuses on this question in a formal way. It provides a first attempt to quantify analytically how the MP convergence is affected by the discretization of the continuous space of dictionary function parameters. 

Our compass to reach this objective is the natural geometry of the continuous dictionary. This dictionary can be seen as a parametric (Riemannian) manifold on which the tools of differential geometry can be applied. This geometrical approach, of increasing interest in the signal processing literature, is inspired by the works \cite{wak05,donoho2005imi} on \emph{Image Appearance Manifolds}, and is also closely linked to manifolds of parametric probability density function associated to the Fisher information metric \cite{amari1982dgc}. Some preliminary hints were also provided in a Riemannian study of generalized correlation of signals with probing functions \cite{watson1998sai}.

The outcome of our study is twofold. On the one hand, we analyze how the rate of convergence of the continuous MP (cMP) is affected by the discretization of the prototype function parameters. We demonstrate that the MP using that discretized dictionary (dMP) converges like a weak continuous MP, i.e. a MP algorithm where the coefficient of the selected atom at each iteration overtakes only a percentage (the \emph{weakness factor}) of the largest atom magnitude. We describe then how this weakness factor decreases as the so-called {\it density radius}\footnote{This density radius represents the maximal distance between any atom of the continuous dictionary and its closest atom in the discretization.} of the discretization increases. This observation is demonstrated experimentally on images and randomly generated 1-D signals. 

On the other hand, to improve the rate of convergence of discrete MP without resorting to a finer but computationally heavier discretization, we propose to exploit a geometric gradient ascent method. This allows to converge to a set of locally optimal continuous parameters, starting from the best set of parameters identified by a coarse but computationally light discrete MP. Each atom of the MP expansion is then defined in two steps. The first step selects the discrete set of parameters that maximizes the inner product between the corresponding dictionary function and the residual signal. The second step implements a (manifold\footnote{In the sense that this gradient ascent evolves on the manifold induced by the intrinsic dictionary geometry.}) gradient ascent method to compute the prototype function parameters that maximize the inner product function over the continuous parameter space. As a main analytical result, we demonstrate that this geometrically optimized discrete MP (gMP) is again equivalent to a continuous MP, but with a weakness factor that is two times closer to unity than for the non-optimized dMP. Our experiments confirm that the proposed gradient ascent procedure significantly increases the rate of convergence of MP, compared to the non-optimized discrete MP. At an equivalent convergence rate, the optimization allows reduction of the discretization density by an order of magnitude, resulting in significant computational gains.

The paper is organized as follows. In Section~2, we introduce the notions of parametric dictionary in the context of signal decomposition in an abstract Hilbert space. This dictionary is then envisioned as a Hilbert manifold, and we describe how its geometrical structure influences its parametrization using the tools of differential geometry. Section~3 surveys the definition of (weak) continuous MP providing a theoretical optimal rate of convergence for further comparisons with other greedy decompositions. A ``discretization autopsy'' of this algorithm is performed in Section~4 and a resulting theorem explaining the dependences of the dMP convergence relatively to this sampling is proved. A simple but illustrative example of a 1-D dictionary, the wavelet (affine) dictionary, is then given. The optimization scheme announced above is developed in Section~5. After a review of gradient ascent optimization evolving on manifolds, the geometrically optimized MP is introduced and its theoretical rate of convergence analyzed in a second theorem. Finally, in Section~6, experiments are performed for 1-D and 2-D signal decompositions using dMP and gMP on various regular discretizations of dictionary parametrizations. We provide links to previous related works in Section~7 and conclude with possible extensions in Section~8. 

\section{Dictionary, Parametrization and Differential Geometry}
\label{sec:dico-param-manifold}

Our object of interest throughout this paper is a general real ``signal'', i.e. a real function $f$ taking value on a measure space $X$. More precisely, we assume $f$ in the set of finite energy signals, i.e. $f\in L^2(X,\ud\mu)=\{u:X\to\Rbb\,:\ \|u\|^2=\int_{X}\,|u(x)|^2\ \ud\mu(x)\ <\ \infty\}$, for a certain integral measure $\ud\mu(x)$. Of course, the natural comparison of two functions $u$ and $v$ in $L^2(X,\ud\mu)$ is realized through the scalar product $\scp{u}{v}_{L^2(X)}=\scp{u}{v}\triangleq\int_{X}\,u(x)\,v(x)\ud\mu(x)$ making $L^2(X,\ud\mu)$ a Hilbert\footnote{Assuming it \emph{complete}, i.e. every Cauchy sequence converges in this space relatively to the norm $\|\cdot\|^2=\scp{\cdot}{\cdot}$.} space where $\|u\|^2=\scp{u}{u}$. 

This very general framework can be specialized to 1-D signal or image decomposition where $X$ is given respectively by $\Rbb$ or $\Rbb^2$, but also to more special spaces like the two dimensional sphere $S^2$ \cite{tosfrovan06} or the hyperboloid \cite{bovdgga05}. In the sequel, we will write simply $L^2(X)=L^2(X,\ud\mu)$.

In the following sections, we will \emph{decompose} $f$ over a highly redundant parametric \emph{dictionary} of real \emph{atoms}. These are obtained from smooth transformations of a real mother function $g\in L^2(X)$ of unit norm. Formally, each atom is a function $g_\lambda(x) = [U(\lambda) g](x) \in L^2(X)$, for a certain isometric operator $U$ parametrized by elements $\lambda\in\Lambda$ and such that $\|g_\lambda\|=\|g\|=1$. The \emph{parametrization} set $\Lambda$ is a continuous space where each $\lambda\in\Lambda$ corresponds to $P$ continuous components $\lambda=\{\lambda^i\}_{0\leq i\leq P-1}$ of different nature. For instance, in the case of 1-D signal or image analysis, $g$ may be transformed by translation, modulation, rotation, or (anisotropic) dilation operations, each associated to one component $\lambda^i$ of $\lambda$. Our dictionary is then the set $\dict(g, U, \Theta) \triangleq \big\{\,g_\lambda(x) = [U(\lambda)g](x): \lambda\in\Theta\,\big\}$, for a certain subset $\Theta \subseteq \Lambda$. In the rest of the paper, we write $\dict(\Theta)=\dict(g,U,\Theta)$, assuming $g$ and $U$ implicitly given by the context. For the case $\Theta=\Lambda$, we write $\mathcal{D}=\dict(\Lambda)$.

We assume that $g$ is twice differentiable over $X$ and that the functions $g_\lambda(x)$ are twice differentiable on each of the $P$ components of $\lambda$. In the following, we write $\partial_i$ for the partial derivative with respect to $\lambda^i$, i.e. $\frac{\partial}{\partial\lambda^i}$, of any element (e.g. $g_\lambda(x)$, $\scp{g_\lambda}{u}$, \ldots) depending on $\lambda$, and $\partial_{ij}=\partial_i\partial_j$. From the smoothness of $U$ and $g$, we have $\partial_{ij}=\partial_{ji}$ on quantities built from these two ingredients. 

Let us now analyze the geometrical structure of $\Lambda$. Rather than an artificial \emph{Euclidean distance} $d_{\mathcal{E}}(\lambda_a,\lambda_b)^2\ \triangleq\ \sum_i(\lambda_a^i-\lambda_b^i)^2$ between $\lambda_a,\lambda_b\in\Lambda$, we use a distance introduced by the dictionary $\mathcal{D}$ itself seen as a $P$-dimensional parametric submanifold of $L^2(X)$ (or a \emph{Hilbert manifold}\,\footnote{This is a special case of Image Appearance Manifold (IAM) defined for instance in \cite{wak05,donoho2005imi}. It is also closely linked to manifolds of parametric probability density function associated to the Fisher information metric \cite{amari1982dgc}.} \cite{lang1972dm}).  The \emph{dictionary distance} $d_{\mathcal{D}}$ is thus the distance in the embedding space $L^2(X)$, i.e. $d_{\mathcal{D}}(\lambda_a,\lambda_b)\ \triangleq\ \|g_{\lambda_a}-g_{\lambda_b}\|$.

From this embedding, we can define an intrinsic distance in $\mathcal{D}$, namely the \emph{geodesic distance}. This later has been used in a similar context in the work of Grimes and Donoho \cite{donoho2005imi} and we follow here their approach. For our two points $\lambda_a,\lambda_b$, assume that we have a smooth curve $\gamma:[0,1]\to\Lambda$ with $\gamma(t)=\big(\gamma^0(t),\cdots,\gamma^{P-1}(t)\big)$, such that $\gamma(0)=\lambda_a$ and $\gamma(1)=\lambda_b$. The length $\mathcal{L}(\gamma)$ of this curve in $\mathcal{D}$ is thus given by
$\mathcal{L}(\gamma) \triangleq \int_0^1\, \|\tfrac{\ud}{\ud t}\,g_{\gamma(t)}\|\,\ud t$, assuming that $g_{\gamma(t)}$ is differentiable\footnote{Another definition of $\mathcal{L}$ exists for non differentiable curve. See for instance \cite{donoho2005imi}.} with respect to $t$.

The \emph{geodesic distance} between $\lambda_a$ and $\lambda_b$ in $\Lambda$ is the length of shortest path between these two points, i.e.
\begin{equation}
\label{eq:geodesic-distance-dico-def}
d_{\mathcal{G}}(\lambda_a,\lambda_b)\ \triangleq\ \inf_{\gamma(\lambda_a\to\lambda_b)}\ \int_0^1 \|\tfrac{\ud}{\ud t}\,g_{\gamma(t)}\|\,\ud t,
\end{equation}
where $\gamma(\lambda_a\to\lambda_b)$ is any differentiable curve $\gamma(t)$ linking $\lambda_a$ to $\lambda_b$ for $t$ equals to 0 and 1 respectively.

We denote by $\gamma_{_{\lambda_a\lambda_b}}$ the optimal \emph{geodesic} curve joining $\lambda_a$ and $\lambda_b$ on the manifold $\mathcal{D}$, i.e. such that $\mathcal{L}(\gamma_{_{\lambda_a\lambda_b}}) = d_{\mathcal{G}}(\lambda_a,\lambda_b)$, and we assume henceforth that it is always possible to define this curve between two points of $\Lambda$. Note that by construction, $d_{\mathcal{G}}(\lambda_a,\lambda_b) = d_{\mathcal{G}}(\lambda_a,\lambda') + d_{\mathcal{G}}(\lambda',\lambda_b)$, for all $\lambda'$ on the curve $\gamma_{_{\lambda_a\lambda_b}}(t)$. 

In the language of differential geometry, the parameter space $\Lambda$ is a Riemannian manifold $\mathcal{M}=(\Lambda,\mathcal{G}_{ij})$ with metric $\mathcal{G}_{ij}(\lambda)=\scp{\partial_i g_\lambda}{\partial_j g_\lambda}$. Indeed, for any differentiable curve $\gamma:t\in[-\delta,\delta] \to \gamma(t)\in\Lambda$ with $\delta>0$ and $\gamma(0)=\lambda$,  we have 
\begin{equation}
\label{eq:pullback-relation-fundation} 
\|\tfrac{\ud}{\ud t}\,g_{\gamma(t)}\big|_{t=0}\|^2\ =\ \dot{\gamma}^i(0)\,\dot{\gamma}^j(0)\,\mathcal{G}_{ij}(\lambda),
\end{equation}
with $\dot{u}(t)=\frac{\ud}{\ud t} u(t)$, and where Einstein's summation convention is used for simplicity\footnote{Namely, a summation in an expression is defined implicitly each time the same index is repeated once as a subscript and once as a superscript, the range of summation being always $[0,P-1]$, so  that for instance the expression $a^i b_i$ reads $\sum_{i=0}^{P-1} a^ib_i$.}. 

The vector $\xi^i=\dot{\gamma}^i(0)$ is by definition a vector in the tangent space $T_\lambda\Lambda$ of $\Lambda$ in $\lambda$. The meaning of relation \eqref{eq:pullback-relation-fundation} is that the metric $\mathcal{G}_{ij}(\lambda)$ allows the definitions of a scalar product and a norm in each $T_\lambda\Lambda$. The norm of a vector $\xi\in T_\lambda\Lambda$ is therefore noted $|\xi|^2 = |\xi|^2_\lambda \triangleq \xi^i\xi^j \mathcal{G}_{ij}(\lambda)$, with the correspondence $\|\tfrac{\ud}{\ud t}\,g_{\gamma(t)}|_{t=0}\| = |\dot{\gamma}|$. For the consistency of further Riemannian geometry developments, we assume that our dictionary $\mathcal{D}$ is \emph{non-degenerate}, i.e. that it induces a positive definite metric $\mathcal{G}_{ij}$ . Appendix~\ref{app:comp-diff-geo} provides additional details.

We conclude this section with the \emph{arc length} (or \emph{curvilinear}) parametrization ``$s$'' \cite{carmo1992rg} of a curve $\gamma(s)$. It is such that $|\gamma'|^2 \triangleq {\gamma'}^i(s)\,{\gamma'}^j(s)\,\mathcal{G}_{ij}(\gamma(s)) = 1$, where $u'(s)=\frac{\ud}{\ud s}u(s)$.  From its definition, the curvilinear parameter $s$ is the one which measures at each point $\gamma(s)$ the length of the segment of curve already travelled on $\gamma$ from $\gamma(0)$. Therefore, in this parametrization, $\lambda_a=\gamma_{_{\lambda_a\lambda_b}}(0)$ and $\lambda_b=\gamma_{_{\lambda_a\lambda_b}}(d_{\mathcal{G}}(\lambda_a,\lambda_b))$. 

\section{Matching Pursuit in Continuous Dictionary}
\label{sec:cMP}

Let us assume that we want to decompose a function $f\in L^2(X)$ into simpler elements (atoms) coming from a dictionary $\dict(\Theta)$, given a possibly uncountable and infinite subset $\Theta\subseteq\Lambda$. Our general aim is thus to find a set of \emph{coefficients} $\{c_m\}$ such that $f(x)$ is equal or well approximated by $f_{\rm app}(x) = \sum_m c_m\,g_{\lambda_m}(x)$ with a finite set of atoms $\{g_{\lambda_m}\}\subset \dict(\Theta)$. 

Formally, for a given \emph{weakness} factor $\alpha\in(0,1]$, a \emph{General Weak$(\alpha)$ Matching Pursuit} decomposition of $f$ \cite{maza93,grivan06}, written ${\rm MP}(\Theta,\alpha)$, in the dictionary $\dict(\Theta)$ is performed through the following \emph{greedy}\footnote{Greedy in the sense that it does not solve a global $\ell_0$ or $\ell_1$ minimization \cite{basispursuit} to find the coefficients $c_m$ of $f_{\rm app}$ above, but works iteratively by solving at each iteration step a local and smaller minimization problem.} algorithm~:\begin{subequations}\label{eq:mp_def}\begin{eqnarray}
&R^0f = f,\ A^0f = 0,\ {\rm(initialization)},\nonumber\\
\label{eq:resid_def}
&R^{m+1}f\ =\ R^{m}f\ -\ \scp{g_{\lambda_{m+1}}}{R^m f}\,g_{\lambda_{m+1}},\ \\
\label{eq:approx_def}
&A^{m+1}f\ =\ A^{m}f\ +\ \scp{g_{\lambda_{m+1}}}{R^m f}\,g_{\lambda_{m+1}},\\
\label{eq:lambda_m_select}
&{\rm with :}\ \scp{g_{\lambda_{m+1}}}{R^m f}^2\ \geq\ \alpha^2\,\mathop{\rm sup}_{\lambda\in\Theta}\,\scp{g_\lambda}{R^m
f}^2.
\end{eqnarray}
\end{subequations}
The quantity $R^{m+1} f$ is the \emph{residual} of $f$ at iteration $m+1$. Since it is orthogonal to atom $g_{\lambda_{m+1}}$, $\|R^{m+1}f\|^2 = \|R^{m}f\|^2 - \scp{g_{\lambda_{m+1}}}{R^mf}^2 \leq \|R^{m}f\|^2$, so that the energy $\|R^mf\|^2$ is non-increasing. The function $A^{m}f$ is the $m$-term \emph{approximation} of $f$ with
$A^{m}f\ =\ \sum_{k=0}^{m-1}\ \scp{g_{\lambda_{k+1}}}{R^k f}\ g_{\lambda_{k+1}}$.

Notice that the \emph{selection rule} \eqref{eq:lambda_m_select} concerns the square of the real scalar product $\scp{g_\lambda}{R^mf}$. Matching Pursuit atom selection is typically defined over the absolute value $|\scp{g_\lambda}{R^mf}|$. However, we prefer this equivalent quadratic formulation first to avoid the abrupt behavior of the absolute value when the scalar product crosses zero, and second for consistency with the quadratic optimization framework to be explained in Section \ref{sec:optim-discr-mp}. Finally, to allow the non-weak case where $\alpha=1$, we assume that a maximizer $g_u\in\dict(\Theta)$ of $\scp{g}{u}^2$ always exists for any $u\in L^2(X)$.
\medskip

If $\Theta$ is uncountable, our general Matching Pursuit algorithm is named \emph{continuous Matching pursuit}. In particular, for $\Theta=\Lambda$, we write ${\rm cMP}(\alpha)={\rm MP}(\Lambda,\alpha)$. The \emph{rate of convergence} (or convergence) of the cMP$(\alpha)$, characterized by the rate of decay of $\|R^m f\|$ with $m$, can be assessed in certain particular cases. For instance, if there exists a Hilbert space $\mathcal{S}\subseteq L^2(X)$ containing $\mathcal{D}=\dict(\Lambda)$ such that
\begin{equation}
\label{eq:beta-space-def}
\beta^2\ =\ \mathop{\rm inf}_{u \in S,\ \|u\|=1}\ \mathop{\rm sup}_{\lambda\in\Lambda}\ \scp{g_\lambda}{u}^2\ >\ 0,
\end{equation}
then the cMP$(\alpha)$ converges inside $\mathcal{S}$. In fact, the convergence is exponential \cite{mall98} since  
$\scp{g_{\lambda_m}}{R^{m-1}f}^2 \geq \alpha^2\beta^2\,\|R^{m-1}f\|^2$ and $\|R^{m}f\|^2 \leq \|R^{m-1}f\|^2 - \alpha^2\beta^2\|R^{m-1}f\|^2 \leq (1-\alpha^2\beta^2)^m\|f\|^2$. We name $\beta=\beta(\mathcal{S},\mathcal{D})$ the \emph{greedy factor} since it charaterizes the MP convergence (greediness).

The existence of the greedy factor $\beta$ is obvious for instance for finite dimensional space \cite{mall98}, i.e. $f\in \Cbb^N$, with finite dictionary (finite number of atoms). 

For a finite dictionary in an infinite dimensional space, as $L^2(X)$, the existence of $\beta$ is not guaranteed over the whole space. However, there exists on the space of functions given by linear combination of dictionary elements, the number of terms being restricted by the dictionary \emph{(cumulative) coherence}~\cite{grivan06}. 

In the case of an infinite dictionary in an infinite dimension space where the greedy factor vanishes, cMP$(\alpha)$ convergence is characterized differently on the subspace of linear combination of countable subsets of dictionary elements. This question is addressed separately in a companion Technical Report \cite{TRLJ0701} to this article. We now consider only the case where a non-zero greedy factor exists to characterize the rate of convergence of MP using continuous and discrete dictionaries.

\section{Discretization effects of Continuous Dictionary}
\label{sec:discr-cont-dico-effects}

The greedy algorithm cMP$(\alpha)$ using the dictionary $\mathcal{D}$ is obviously numerically unachievable because of the intrinsic continuity of its main ingredient, namely the parameter space $\Lambda$. Any computer implementation needs at least to discretize the parametrization of the dictionary, more or less densely, leading to a countable set $\Lambda_{\rm d} \subset \Lambda$. This new parameter space leads naturally to the definition of a countable subdictionary $\mathcal{D}_{\rm d}=\dict(\Lambda_\ud)$. Henceforth, elements of $\Lambda_\ud$ are labelled with roman letters, e.g. $k$, to distinguish them from the continuous greek-labelized elements of $\Lambda$, e.g $\lambda$.

For a weakness factor $\alpha\in(0,1]$, the \emph{discrete Weak$(\alpha_{\rm d})$ Matching Pursuit} algorithm, or dMP$(\alpha)$, of a function $f\in L^2(X)$ over $\mathcal{D}_{\rm d}$ is naturally defined as ${\rm dMP}(\alpha) = {\rm MP}(\Lambda_\ud,\alpha)$.
The replacement of $\Lambda$ by $\Lambda_\ud$ in the MP algorithm \eqref{eq:mp_def} leads obviously to the following question that we address in the next section. 

\begin{question}
\label{quest:discr-effect}
How does the MP rate of convergence evolve when the parametrization of a dictionary is discretized and what are the quantities that control (or bound) this evolution~?
\end{question}

\subsection{Discretization Autopsy} 

By working with $\mathcal{D}_{\rm d}$ instead of $\mathcal{D}$, the atoms selected at each iteration of dMP$(\alpha)$ are of course less optimal than those available in the continuous framework. Answering Question \ref{quest:discr-effect} requires a quantitative measure of the induced loss in the MP coefficients. More concretely, defining the \emph{score function} $S_u(\lambda)=\scp{g_\lambda}{u}^2$ for some $u\in L^2(X)$, we must analyze the difference between a maximum of $S_u$ computed over $\Lambda$ and that obtained from $\Lambda_{\rm d}$. This function $u$ will be next identified with the residue of dMP$(\alpha)$ at any iteration to characterize the global change in convergence.

We propose to found our analysis on the geometric tools described in Section \ref{sec:dico-param-manifold}. 

\begin{definition}
The value $S_u(\lambda_a)$ is critical in the direction of $\lambda_b$ if, given the geodesic $\gamma=\gamma_{_{\lambda_a\lambda_b}}$ in the manifold $\mathcal{M}=(\Lambda,\mathcal{G}_{ij})$, $\frac{\ud}{\ud s}S_u(\gamma(s))|_{s=0} = 0$, where $\gamma(0)=\lambda_a$.
\end{definition}

Notice that if $S_u(\lambda_a)$ is critical in the direction of $\lambda_b$, ${\gamma'}^i(0)\,\partial_i S_u(\lambda_a) = 0$. An \emph{umbilical} point for which $\partial_i S_u(\lambda_a) = 0$ for all $i$, is obviously critical in any direction. An umbilical point corresponds geometrically either to maxima, minima or saddlepoints of $S_u$ relatively to $\Lambda$. 

\begin{proposition}
\label{prop:quad-approx-best-scp}
Given $u\in L^2(X)$, if $S_u(\lambda_a)$ is critical in the direction of $\lambda_b$ for $\lambda_a,\lambda_b\in\Lambda$, then for some $r\in (0, d_{\mathcal{G}}(\lambda_a,\lambda_b))$,
\begin{equation}
\label{eq:quad-approx-best-scp}
|S_u(\lambda_a) - S_u(\lambda_b)| \leq \|u\|^2\,d_{\mathcal{G}}(\lambda_a,\lambda_b)^2\,\big(1 + \|\tfrac{\ud^2g_{\gamma}}{\ud s^2}\big\vert_{s=r}\|\big),
\end{equation}
where $\gamma(s)=\gamma_{_{\lambda_a\lambda_b}}(s)$ is the geodesic in $\mathcal{M}$ linking $\lambda_a$ to $\lambda_b$.
\end{proposition}

\begin{proof}
Let us define the twice differentiable function $\psi(s) \triangleq S_u(\gamma(s))$ on $s\in [0,\eta]$, with $\eta\triangleq d_{\mathcal{G}}(\lambda_a,\lambda_b)$. A second order Taylor development of $\psi$ gives, for a certain $r\in(0,s)$,
$\psi(s) = \psi(0) + s\,\psi'(0) + \tinv{2}s^2\,\psi''(r)$. Since $\psi'(0)={\gamma'}^i(0)\,\partial_iS_u(\lambda_a)=0$ by hypothesis, we have in $s=\eta$, $|\psi(0)-\psi(\eta)| = |S_u(\lambda_a)-S_u(\lambda_b)| \leq \tinv{2}\,\eta^2\,|\psi''(r)|$. However, on any $s$, $|\psi''(s)| = 2\,|\bscp{\tfrac{\ud}{\ud s}g_{\gamma(s)}}{u}^2 + \scp{g_{\gamma(s)}}{u}\,\bscp{\tfrac{\ud^2}{\ud s^2}g_{\gamma(s)}}{u}| \leq 2\,(\|\tfrac{\ud}{\ud s}g_{\gamma(s)}\|^2 + \|\tfrac{\ud^2}{\ud s^2}g_{\gamma(s)}\|)\,\|u\|^2$, using the Cauchy-Schwarz (CS) inequality in $L^2(X)$ in the last equation. The result follows from the fact that $\|\tfrac{\ud}{\ud s}g_{\gamma(s)}\|=1$.
\end{proof}

The previous Lemma is particularly important since it bounds the loss in coefficient value when we decide to choose $S_u(\lambda_b)$ instead of the optimal $S_u(\lambda_a)$ in function of the geodesic distance $d_{\mathcal{G}}(\lambda_a,\lambda_b)$ between the two parameters. To obtain a more satisfactory control of this difference, we need however a new property of the dictionary. 

We start by defining the \emph{principal curvature} in the point $\lambda\in\Lambda$ as 
\begin{equation}
\label{eq:principal-curv-def}
\mathcal{K}_\lambda\ \triangleq\ \sup_{\xi\,:\ |\xi|=1}\,\|\tfrac{\ud^2}{\ud s^2}\,g_{\gamma_\xi(s)}\big\vert_{s=0}\|,
\end{equation}
where $\gamma_\xi$ is the unique geodesic in $\mathcal{M}$ starting from $\lambda=\gamma_\xi(0)$ and with $\gamma'_\xi(0)=\xi$, for a direction $\xi$ of unit norm in $T_\lambda\Lambda$. 

\begin{definition}
\label{def:bound-curve-dico}
The \emph{condition number} of a dictionary $\mathcal{D}$ is the number $\mathcal{K}^{-1}$ obtained from 
\begin{equation}
\label{eq:bound-curve-dico}
\mathcal{K}\ \triangleq\ \sup_{\lambda\in\Lambda}\,\mathcal{K}_\lambda.
\end{equation}
If $\mathcal{K}$ does not exist (not bounded $\mathcal{K}_\lambda$), by extension, $\mathcal{D}$ is said to be of zero condition number.
\end{definition}

The notion of condition number has been introduced by Niyogi et al. \cite{niyogi2004fhs} to bound the local curvature of an embedded manifold\footnote{In their work, the condition number, named there $\tau^{-1}$, of a manifold $\mathcal{M}'$ measures the maximal ``thickness'' $\tau$ of the \emph{normal bundle}, the union of all the orthogonal complement of every tangent plane at every point of the manifold.} in its ambient space, and to characterize its self-avoidance. Essentially, it is the inverse of the maximum radius of a sphere that, when placed tangent to the manifold at any point, intersects the manifold only at that point \cite{davenport2007sfc,wak06manif}. Our quantity $\mathcal{K}^{-1}$ is then by construction a similar notion for the dictionary $\mathcal{D}$ seen as a manifold in $L^2(X)$. However, it does not actually prevent manifold self-crossing on large distance due to the locality of our differential analysis\footnote{A careful study of local self-avoidance of well-conditioned dictionary would have to be considered but this is beyond the scope of this paper.}.

\begin{proposition}
\label{prop:bounds-dico-curv} 
For a dictionary $\mathcal{D}=\dict(\Lambda)$, 
\begin{equation}
\label{eq:unit-lower-bound-dico-curv}
1\ \leq\ \mathcal{K}\ \leq\ \sup_{\lambda\in\Lambda}\ \left[\,\bscp{\partial_{ij}\,g_{\lambda}}{\partial_{kl}\,g_{\lambda}}\,\mathcal{G}^{ik}\,\mathcal{G}^{jl}\,\right]^{\tinv{2}},
\end{equation}
where $\mathcal{G}^{ij}=\mathcal{G}^{ij}(\lambda)$ is the inverse\footnote{Using Einstein convention, this means $\mathcal{G}^{ik}\mathcal{G}_{kj}=\mathcal{G}_{jk}\mathcal{G}^{ki}=\delta^i_j$, for the Kronecker's symbol $\delta^i_j=\delta_{ij}=\delta^{ij}=1$ if $i=j$ and 0 if $i\neq j$.} of $\mathcal{G}_{ij}$. 
\end{proposition}
The proof is given in Appendix \ref{app:proof-prop-unit-lower-bound-dico-curv} since it uses some elements of differential geometry not essential in the core of this paper. The interested reader will find also there a slightly lower bound than the bound presented in \eqref{eq:unit-lower-bound-dico-curv}, exploiting covariant derivatives, Laplace-Beltrami operator and scalar curvature of $\mathcal{M}$~\cite{carmo1992rg}. We can state now the following corollary of Proposition \ref{prop:quad-approx-best-scp}. 

\begin{corollary}
\label{coro:bound-scp-diff-for-bound-curv-dico}
In the conditions of Proposition \ref{prop:quad-approx-best-scp}, if $\mathcal{D}$ has a non-zero condition number $\mathcal{K}^{-1}$, then 
\begin{equation}
\label{eq:quad-approx-best-scp-cond-number}
|S_u(\lambda_a) - S_u(\lambda_b)| \quad \leq \quad \|u\|^2\,d_{\mathcal{G}}(\lambda_a,\lambda_b)^2\,\big(1 + \mathcal{K}\big).
\end{equation}
\end{corollary}

Therefore, in the dMP$(\alpha)$ decomposition of $f$ based on $\mathcal{D}_{\rm d}$, even if at each iteration the exact position of the continuous optimal atom of $\mathcal{D}$ is not known, we are now able to estimate the convergence rate of this MP provided we introduce a new quantity characterizing the set $\Lambda_{\rm d}$. 

\begin{definition}
\label{def:covering-def-density-radius}
The \emph{density radius} $\rho_{\rm d}$ of a countable parameter space $\Lambda_{\rm d}\subset\Lambda$ is the value
\begin{equation}
\label{eq:dens-radius-def}
\rho_{d}\ =\ \sup_{\lambda\in\Lambda}\,\inf_{k\in\Lambda_{\rm d}}\ d_{\mathcal{G}}(\lambda,\,k).
\end{equation}
We say that $\Lambda_{\rm d}$ \emph{covers} $\Lambda$ with a radius $\rho_{\rm d}$.
\end{definition}
This radius characterizes the density of $\Lambda_{\rm d}$ inside $\Lambda$. Given any $\lambda$ in $\Lambda$, one is guaranteed that there exists an element $k$ of $\Lambda_{\rm d}$ close to $\lambda$, i.e. within a geodesic distance $\rho_{\rm d}$.

\begin{theorem}
\label{prop:disc_dico-density-requ-conver-rate}
Given a Hilbert space $\mathcal{S}\subseteq L^2(X)$ with a non zero greedy factor $\beta$, and a dictionary $\mathcal{D}=\dict(\Lambda)\subset S$ of non-zero condition number $\mathcal{K}^{-1}$, if $\Lambda_{\rm d}$ covers $\Lambda$ with radius $\rho_\ud$, and if $\rho_{\rm d} < \beta/\sqrt{1+\mathcal{K}}$, then, for functions belonging to $\mathcal{S}$, a dMP$(\alpha)$ algorithm using $\mathcal{D}_{\rm d}=\dict(\Lambda_\ud)$ is bounded by the exponential convergence rate of a cMP$(\alpha')$ using $\mathcal{D}$ with a weakness parameter given by $\alpha' = \alpha\big(1-\beta^{-2}\,\rho_{\rm d}^2(1+\mathcal{K})\big)^{1/2} <\alpha$.
\end{theorem}

\begin{proof} Notice first that since $f\in\mathcal{S}$ and $\mathcal{D}_\ud\subset\mathcal{D}\subset\mathcal{S}$, $R^m f\in\mathcal{S}$ for all iteration $m$ of dMP. Let us take the $(m+1)^{\rm th}$ step of dMP$(\alpha)$ and write $u=R^m f$.  We have of course $\|R^{m+1} f\|^2=\|u\|^2-S_u(k_{m+1})$, where $k_{m+1}$ is the atom obtained from the selection rule \eqref{eq:lambda_m_select}, i.e. $S_u(k_{m+1}) \geq \alpha^2\,\sup_{k\in\Lambda_{\rm d}}\,S_u(k)$. 

Denote by $g_{\tilde{\lambda}}$ the atom of $\mathcal{D}$ that best represents $R^m f$, i.e. $S_u(\tilde{\lambda}) = \sup_{\lambda\in\Lambda}S_u(\lambda)$. If $\tilde{k}$ is the closest element of $\tilde{\lambda}$ in $\Lambda_{\rm d}$, we have $d_{\mathcal{G}}(\tilde{\lambda},\tilde{k})\leq\rho_{\rm d}$ from the covering property of $\Lambda_{\rm d}$, and the Proposition \ref{prop:quad-approx-best-scp} tells us that, with $u=R^m f$, $|S_u(\tilde{k})-S_u(\tilde{\lambda})| \leq \rho_{\rm d}^2\,(1+\mathcal{K})\,\|u\|^2$, since $\partial_i S_u(\tilde{\lambda})=0$ for all $i$.

Therefore, $S_u(\tilde{k}) \geq S_u(\tilde{\lambda}) - \rho_{\rm d}^2\,(1+\mathcal{K})\,\|u\|^2 \geq \beta^2\,\|u\|^2 - \rho_{\rm d}^2\,(1+\mathcal{K})\,\|u\|^2$, and $S_u(\tilde{k}) \geq \beta^2\,\big(1-\beta^{-2}\,\rho_{\rm d}^2(1+\mathcal{K})\big)\,\|R^m f\|^2$, this last quantity being positive from the density requirement, i.e. $\rho_{\rm d} < \beta/\sqrt{1+\mathcal{K}}$.

In consequence, $S_u(k_{m+1}) \geq \alpha^2\,\sup_{k\in\Lambda_{\rm d}}\,S_u(k) \geq \alpha^2\,S_u(\tilde{k})$, implying $\|R^{m+1} f\|^2\ =\ \|u\|^2\ -\ S_u(k_{m+1})\ \leq\ \|u\|^2\ -\ \alpha^2\,S_u(\tilde{k})\ \leq\ \|u\|^2\,(1 - \alpha'^2\beta^2)$, for $\alpha' \triangleq \alpha\big(1-\beta^{-2}\,\rho_{\rm d}^2(1+\mathcal{K})\big)^{1/2}$. So, $\|R^{m+1} f\|\leq(1-\alpha'^2\beta^2)^{(m+1)/2}\|f\|$, which is the exponential convergence rate of the Weak$(\alpha)$ Matching Pursuit in $\mathcal{D}$ when $\beta$ exists \cite{mall98,grivan06}.
\end{proof}

The previous proposition has an interesting interpretation : a weak Matching Pursuit decomposition in a discrete dictionary corresponds, in terms of rate of convergence, to a weaker Matching Pursuit in the continuous dictionary from which the discrete one is extracted. 

About the hypotheses of the proposition, notice first that the existence of a greedy factor inside $\mathcal{S}$ concerns the continuous dictionary $\mathcal{D}$ and not the discrete one $\mathcal{D}_\ud$. Consequently, this condition is certainly easier to fulfill from the high redundancy of $\mathcal{D}$. Second, the \emph{density requirement}, $\rho_{\rm d} < \beta/\sqrt{1+\mathcal{K}}$, is just sufficient since the Proposition \ref{prop:quad-approx-best-scp} does not state that it achieves the best bound for the control of $|S_u(\lambda_a)-S_u(\lambda_b)|$ when $\lambda_a$ is critical. It is interesting to note that this inequality relates $\rho_{\rm d}$, a quantity that characterizes the discretization $\Lambda_{\rm d}$ in $\Lambda$, to $\beta$ and $\mathcal{K}$, which depend only on the dictionary. In particular, $\beta$ represents the density of $\mathcal{D}$ inside $\mathcal{S}\subset L^2(X)$, and $\mathcal{K}$ depends on the shape of the atoms through the curvature of the dictionary.

Finally note that as $\beta<1$ (from definition \eqref{eq:beta-space-def}) and $\mathcal{K}>1$ (Prop. \ref{prop:bounds-dico-curv}), the density radius must at least satisfy $\rho_{\rm d} < \inv{\sqrt{2}}$ to guarantee that our analysis is valid. 

\subsection{A Simple Example of Discretization} 
\label{subsec:example}

Let us work on the \emph{line} with $L^2(X)=L^2(\Rbb,\ud t)$, and check if the hypothesis of the previous theorem can be assessed in the simple case of an \emph{affine} (wavelet-like) dictionary. 

We select a symmetric and real mother function $g\in L^2(\Rbb)$ well localized around the origin, e.g. a Gaussian or a Mexican Hat, normalized such that $\|g\|=1$. The parameter set $\Lambda$ is related to the \emph{affine group}, the group of translations and dilations $G_{\rm aff}$. We identify $\lambda=(\lambda^0=b,\lambda^1=a)$, where $b\in\Rbb$ and $a>0$ are the translation and the dilation parameters respectively. The dictionary $\mathcal{D}$ is defined from the atoms $g_\lambda(t) = [U(\lambda)g](t) = a^{-1/2}\,g\big((t-b)/a\big)$, with $\|g_\lambda\|=1$ for all $\lambda\in\Lambda$. Our atoms are nothing but the wavelets of a Continuous Wavelet Transform if $g$ is admissible \cite{daub92}, and $U$ is actually the representation of the affine group on $L^2(\Rbb)$ \cite{alanga99}. 

In the technical report \cite{TRLJ0701}, we prove that the associated metric is given by $\mathcal{G}_{ij}(\lambda)=a^{-2}\,W$, where $W$ is a constant $2\times 2$ diagonal matrix depending only of the mother function $g$ and its first and second derivatives. Since $\mathcal{G}^{ij}(\lambda)=a^{2}\,W^{-1}$, $\mathcal{K}$ can be bounded by a constant also associated to $g$ and its first and second order time derivatives.

Finally, given the $\tau$-adic parameter discretization
$$
\Lambda_{\rm d}=\{k_{jn}=(b_{jn},a_j)=(n\,b_0\,\tau^j,a_0\tau^j):\ j,n\in\Zbb\},
$$
with $\tau>1$ and $a_0,b_0>0$, the density radius $\rho_{\rm d}$ of $\Lambda_{\rm d}$ is shown to be bounded by $\rho_{\rm d}\leq C a_0^{-1} b_0 + D\ln\tau$, with $C$ and $D$ depending only of the norms of $g$ and its first derivative. 

This bound has two interesting properties. First, as for the grid $\Lambda_\ud$, it is invariant under the change $(b_0,a_0)\to (2b_0, 2a_0)$. Second, it is multiplied by $2^n$ if we realize a ``zoom'' of factor $2^n$ in our $\tau$-adic grid, in other words, if $(b_0,\tau) \to (2^n\,b_0, \tau^{2^n})$. By the same argument, the true density radius has also to respect these rules. Therefore, we conjecture that $\rho_{\rm d} = C' a_0^{-1} b_0 + D'\,\ln\tau$, for two particular (non computed) positive constants $C'$ and $D'$.

Unfortunately, even for this simple affine dictionary, the existence of $\beta=\beta(\mathcal{S},\mathcal{D})$ is non trivial to prove. However, if the greedy factor exists, the control of $\tau$, $a_0$ and $b_0$ over $\rho_{\rm d}$ tells us that it is possible to satisfy the density requirement for convenient values of these parameters. 

\section{Optimization of Discrete Matching Pursuits}
\label{sec:optim-discr-mp}

The previous section has shown that under a few assumptions a dMP is equivalent, in terms of rate of convergence, to a weaker cMP in the continuous dictionary from which the discrete one has been sampled. 
\begin{question}
\label{quest:improve-discr-effect}
Can we improve the rate of convergence of a dMP, not with an obvious increasing of the dictionary sampling, but by taking advantage of the dictionary geometry~?
\end{question}

Our approach is to introduce an optimization of the discrete dMP scheme. In short, at each iteration, we propose to use the atoms of $\mathcal{D}_\ud$ as the seeds of an iterative optimization, such as the basic \emph{gradient descent/ascent}, respecting the geometry of the manifold $\mathcal{M}=(\Lambda,\mathcal{G}_{ij})$. 

Under the same density hypothesis of Theorem \ref{prop:disc_dico-density-requ-conver-rate}, we show that in the worst case and if the number of optimization steps is large enough, an optimized discrete MP is again equivalent to a continuous dMP, but with a weakness factor two times closer to unity than for the non-optimized discrete MP. 

In this section, we first introduce the basic gradient descent/ascent on a manifold. Next, we show how this optimization can be introduced in the Matching Pursuit scheme to defined the geometrically optimized MP (gMP). Finally, the rate of convergence of this method is analyzed. 

\subsection{Gradient Ascent on Riemannian Manifolds}
\label{sec:optim-manifold}

Given a function $u\in L^2(X)$ and $S_u(\lambda)=\scp{g_\lambda}{u}^2$, we wish to find the parameter that maximizes $S_u$, i.e. 
\begin{equation*}
\lambda_*\ =\ \displaystyle\mathop{\rm arg\,max}_{\lambda\in\Lambda}\ S_u(\lambda)
\eqno{\bf (P.1)}
\end{equation*}

Equivalently, by introducing $h_{u,\lambda}=\scp{g_{\lambda}}{u}\,g_\lambda$, we can decide to find $\lambda_*$ by the minimization
\begin{equation*}
\lambda_*\ =\ \displaystyle\mathop{\rm arg\,min}_{\lambda\in\Lambda}\ \|u - h_{u,\lambda}\|^2.
\eqno{\bf (P.2)}
\end{equation*}
If we are not afraid to get stuck on local maxima (P.1) or minima (P.2) of these two non-convex problems, we can solve them
by using well known optimization techniques such as gradient descent/ascent, or Newton or Newton-Gauss optimizations.

We present here a basic gradient ascent of the Problem (P.1) that respect the geometry of $\mathcal{M}=(\Lambda,\mathcal{G}_{ij})$ \cite{ferreira1998sar}. This method increases iteratively the value of $S_u$ by following a path in $\Lambda$, composed of geodesic segments, driven by the gradient of $S_u$.

Given a sequence of step size $t_r>0$, the gradient ascent of $S_u$ starting from $\lambda_0\in\Lambda$ is defined by the following induction \cite{gabay1982mdf} :
$$
\phi_0(\lambda_0)\ =\ \lambda,\quad \phi_{r+1}(\lambda_0)\ =\ \gamma\big(t_r,\ \phi_{r}(\lambda_0),\ \xi_r(\lambda_0)\,\big),
$$
where $\xi_r(\lambda_0) = |\nabla S_u(\phi_r(\lambda_0))|^{-1}\,\nabla S_u(\phi_r(\lambda_0))$ is the \emph{gradient direction} obtained from the gradient $\nabla^i S_u$ $ = \mathcal{G}^{ij}\,\partial_j S_u$, and $\gamma(s,\lambda_0,\xi_0)$ is the geodesic starting at $\lambda_0=\gamma(0,\lambda_0,\xi_0)$ with the unit velocity $\xi_0=\tfrac{\partial}{\partial s}\gamma(0,\lambda_0,\xi_0)$. Notice that $\nabla^i$ is the natural notion of gradient on a Riemannian manifold. Indeed, as for the Euclidean case, with $\nabla^i h \triangleq \mathcal{G}^{ij}\,\partial_{j} h$ for $h\in L^2(X)$, given $w\in T_\lambda\Lambda$, the directional derivative $D_w h$ is equivalent to $D_{w} h(\lambda) \triangleq w^i\partial_i h(\lambda) = \scp{\nabla h}{w}_{\lambda} \triangleq  w^i\,\nabla^j h(\lambda)\,\mathcal{G}_{ij}(\lambda)$, since $\mathcal{G}^{ik}\,\mathcal{G}_{kj}=\delta^{i}_{j}$.

Practically, in our gradient ascent, we use the linear first order approximation of $\gamma$, i.e.
\begin{equation}
\label{eq:grad-asc-optim-func-def}
\phi_{r+1}(\lambda)\ =\ \phi_{r}(\lambda)\ +\ t_r\,\xi_r(\lambda),
\end{equation}
valid for small value of $t_r$ (error in $O(t_r^2)$). This is actually an optimization method since $\partial_i S_u(\phi_{r}(\lambda))\,\xi^i_r = |\partial S_u(\phi_{r}(\lambda))| > 0$ and $S_u(\phi_{r+1}(\lambda)) = S_u(\phi_{r}(\lambda)) + t_r |\partial S_u(\phi_{r}(\lambda))| + O(t_r^2) \geq S_u(\phi_{r}(\lambda))$, for a convenient step size $t_r>0$. At each step of this gradient ascent, the value $t_r$ is chosen so that $S_u$ is increased. This can be done for instance by a \emph{line search} algorithm \cite{bova04}. From the positive definiteness of $\mathcal{G}_{ij}$ and $\mathcal{G}^{ij}$, a fixed point $\phi_{r+1}(\lambda)=\phi_{r}(\lambda)$ is reached if $\nabla^iS_u(\phi_r(\lambda))=\partial_i S_u(\phi_r(\lambda)) =0$ for all $i$.

More sophisticated algorithms such as Newton or Newton-Gauss can be developed to solve the Problem (P.2) on a Riemannian manifolds \cite{mahony2002gnm, gabay1982mdf} even if, unlike to the flat case, a direct definition of the Hessian does not exist on differentiable manifolds. However, we will not use them here as our aim is to prove that a dMP driven by the very basic optimization above provides already a better rate of convergence than the non-optimized dMP.

\subsection{Optimized Discrete Matching Pursuit Algorithm}
\label{subsec:optim-discr-MP-algo}

Let us optimize each step of a discrete MP using the gradient ascent of the previous section.

\paragraph{Definition } Given sequence of positive integers $\kappa_m$ and a weakness factor $0<\alpha\leq 1$, the geometrically optimized discrete matching pursuit (gMP$(\alpha)$) is defined by
\begin{subequations}
\label{eq:mp_optim_disc_def}
\begin{eqnarray}
\label{eq:mp_init_opt}
&R^0 f\ =\ f\quad \textrm{(initialization)},\\
\label{eq:resid_opt_def}
&R^{m+1} f\ =\ R^{m} f\ -\ \scp{g_{\nu_{m+1}}}{R^m f}\,g_{\nu_{m+1}},
\end{eqnarray}
\begin{equation}\label{eq:atom_optim_select}
\scp{g_{\nu_{m+1}}}{R^m f}^2\ \geq\ \alpha^2\,\sup_{k\in\Lambda_{\rm d}}\,\scp{g_{\phi_{\kappa_m}(k)}}{R^m f}^2.
\end{equation}
\end{subequations}
Notice that the best atom $g_{\nu_{m+1}}$ is selected in the set $\Phi_m \triangleq \{g_{\phi_{\kappa_m}(k)}:k\in\Lambda_{\rm d}\}\subset\mathcal{D}$. Elements of $\Phi_m$ are determined by applying the \emph{optimization function} $\phi_r:\Lambda_{\rm d}\to\Lambda$ of our gradient ascent defined in \eqref{eq:grad-asc-optim-func-def} on elements of $\Lambda_{\rm d}$. In consequence, $\Phi_m$ depends on $R^m f$ and is thus different at each iteration $m$. 

\paragraph{Rate of convergence } The following theorem characterizes the rate of convergence of the optimized Matching Pursuit defined in \eqref{eq:mp_optim_disc_def}. 

\begin{theorem}
 \label{prop:optim-dmp-eq-weak-mp}
Given the notations and the conditions of Theorem \ref{prop:disc_dico-density-requ-conver-rate}, there exists a sequence of positive integers $\kappa_m$ such that, the gMP$(\alpha)$ decomposition of functions in $\mathcal{S}\subset L^2(X)$ optimized $\kappa_m$ steps at each iteration $m$, is bounded by the same rate of convergence as a cMP$(\alpha'')$ using the corresponding continuous dictionary $\mathcal{D}$ with $\alpha''=\alpha(1-\inv{2}\,\beta^{-2}\,\rho_{\rm d}\,(1+\mathcal{K}))^{1/2} \leq \alpha$. 
\end{theorem}

In other words, for $\alpha=1$, a gMP is equivalent to a cMP with a weakness factor two times closer to unity than the one reached by a dMP in the same conditions. Before proving this result, let us introduce some new lemmata. 

\begin{lemma}
\label{prop:distance-between-extrema}
Given a function $u\in L^2(X)$ and a dictionary $\mathcal{D}$ of non-zero condition number $\mathcal{K}^{-1}$, if $\lambda_a$ is critical in the direction of $\lambda_b$, and if $\lambda_b$ is critical in the direction of $\lambda_a$, i.e. ${\gamma'}^i(0)\,\partial_iS_u(\lambda_a)={\gamma'}^i(d)\,\partial_iS_u(\lambda_b)=0$ for $\gamma=\gamma_{_{\lambda_a\lambda_b}}$ the geodesic joining $\lambda_a$ and $\lambda_b$ and $d=d_{\mathcal{G}}(\lambda_a,\lambda_b)$, then 
\begin{equation}
 |S_u(\lambda_a)-S_u(\lambda_b)|\ \leq\ \tinv{2}\,\|u\|^2\,d_{\mathcal{G}}(\lambda_a,\lambda_b)^2\,(1+\mathcal{K}).
\end{equation}
\end{lemma}
\begin{proof} Without loss of generality, assume that $S_u(\lambda_a)\geq S_u(\lambda_b)$. If this is not the case, we can switch the labels $a$ and $b$. Let us define $\lambda(\theta)=\gamma(\theta d)$ with $\theta\in [0,1]$ on the geodesic $\gamma=\gamma_{_{\lambda_a\lambda_b}}$. We have $\lambda_a=\lambda(0)$ and $\lambda_b=\lambda(1)$. Using the Corollary \ref{coro:bound-scp-diff-for-bound-curv-dico}, the two following inequalities hold~: $S_u(\lambda(\theta)) \geq S_u(\lambda_a) - \|u\|^2\,d_{\mathcal{G}}(\lambda(\theta), \lambda_a)^2\,(1+\mathcal{K})$ and $S_u(\lambda(\theta)) \leq S_u(\lambda_b) + \|u\|^2\,d_{\mathcal{G}}(\lambda(\theta), \lambda_b)^2\,(1+\mathcal{K})$.

Therefore, since by definition of $\lambda(\theta)$, $d_{\mathcal{G}}(\lambda(\theta), \lambda_a)=\theta d$ and $d_{\mathcal{G}}(\lambda(\theta), \lambda_b)=(1-\theta)d$, we find $S_u(\lambda_a)-S_u(\lambda_b) \leq  \|u\|^2\,\big(\theta^2+(\theta-1)^2\big)\,d^2\,(1+\mathcal{K})$ for all $\theta\in [0,1]$. Taking the minimum over all $\theta$, we obtain finally $S_u(\lambda_a)-S_u(\lambda_b)\ \leq\ \tinv{2}\,\|u\|^2\,d_{\mathcal{G}}(\lambda_a,\lambda_b)^2\,(1+\mathcal{K})$.
\end{proof} 

In other words, the critical nature of $\lambda_a$ and $\lambda_b$ divides by two the bound on the decreasing of $S_u$ between them compared to the situation where only one of these points is critical. 

\begin{lemma}
\label{prop:distance-between-maxima-one-orbit-point}
Given a function $u\in L^2(X)$, assume that $S_u(\lambda)$ has a global maximum at $\lambda_M$, i.e. $\partial_iS_u(\lambda_M)=0$ for all $i$, and write $\mathcal{T}_k=\{\phi_r(k):r\in\Nbb\}$ the trajectory of the gradient ascent described in \eqref{eq:grad-asc-optim-func-def} starting from a point $k\in\Lambda_\ud$. There exists a $\lambda'\in\mathcal{T}_k$ that can be reached in a finite number of optimization steps, such that
\begin{equation}
\label{eq:max-futur-k-rel}
S_u(\lambda_M) - S_u(\lambda')\ \leq\ \tinv{2}\,\|u\|^2\ d_{\mathcal{G}}(\lambda_M,k)^2\,(1+\mathcal{K}).
\end{equation}
\end{lemma}

For the sake of clarity, the proof of this technical Lemma is placed in Appendix \ref{app:proof-distance-between-maxima-one-orbit-point}. The main idea is to find a point in the trajectory $\mathcal{T}_k$ that is closer to $\lambda_M$ than $k$, and that is also critical in the direction of $\lambda_M$ so that Lemma \ref{prop:distance-between-extrema} can be applied. Let us now enter in the proof of the previous proposition.

\begin{proof}[Proof of Theorem \ref{prop:optim-dmp-eq-weak-mp}] In our gMP$(\alpha)$ decomposition of a function $f\in \mathcal{S}\subset L^2(X)$ defined before, given the iteration $m+1$ where $u=R^m f$ is analyzed, denote by $\tilde{\lambda}$ the parameter of the atom in $\mathcal{D}$ maximizing $S_u$, i.e. $S_u(\tilde{\lambda})=\sup_{\lambda\in\Lambda}S_u(\lambda)$.

If $\tilde{k}$ is the closest element of $\Lambda_{\rm d}$ to $\tilde{\lambda}$, from the covering property of $\Lambda_{\rm d}$ we have $d_{\mathcal{G}}(\tilde{\lambda},\tilde{k})\leq \rho_{\rm d}$, and the Lemma \ref{prop:distance-between-maxima-one-orbit-point} tells us that there exists a finite number of optimization steps $\kappa_m$ such that $S_u(\phi_{\kappa_m}(\tilde{k})) \geq S_u(\tilde{\lambda}) - \tinv{2}\,\rho_{\rm d}\,\|u\|^2\,(1+\mathcal{K}) \geq \beta^2\,\big(1-\tinv{2}\beta^{-2}\,\rho_{\rm d}^2\,(1+\mathcal{K})\big)\,\|u\|^2$, where the last term is positive from the density requirement $\rho_{\rm d}< \beta/\sqrt{1+\mathcal{K}}$.

Therefore, from the selection rule \eqref{eq:atom_optim_select}, $S_u(\nu_{m+1}) \geq \alpha^2\,S_u(\phi_{\kappa_m}(\tilde{k}))$. We have thus $\|R^{m+1} f\|^2 = \|u\|^2 - S_u(\nu_{m+1}) \leq \|u\|^2 - \alpha^2\,S_u(\phi_{\kappa_m}(\tilde{k})) \leq \|u\|^2\,(1-\alpha''^2\beta^2)$, with $\alpha''\triangleq\alpha\big(1-\tinv{2}\beta^{-2}\,\rho_{\rm d}^2\,(1+\mathcal{K})\big)^{1/2}$. So, $\|R^{m+1} f\|\leq(1-\alpha''^2\beta^2)^{(m+1)/2}\|f\|$ which is also the exponential convergence rate of the cMP$(\alpha'')$ in $\mathcal{D}$ when $\beta$ exists.
\end{proof}

In Theorem \ref{prop:optim-dmp-eq-weak-mp}, even if the sequence of optimization steps $\kappa_m$ is proved to exist, it is actually unknown. One practical way to overcome this problem is to observe how the ratio $\frac{|\nabla S_u|}{S_u}$ decreases at each optimization steps, and to stop the procedure once this value falls below a predefined threshold. This follows from the idea that the closer to a local maximum $S_u(\phi_r(k))$ is, the smaller must be the optimization step. As it is often the case in optimization problems, an upper bound on the number of optimization steps can be fixed jointly to this threshold test.

\section{Experiments}
\label{sec:experiments}

In this section, dMP and gMP decompositions of 1-D and 2-D signals are studied experimentally in different situations. These will imply different classes of signals and different discretization of parametrization of various densities. 

Prior to these experiments, some remarks have to be made about dMP and gMP implementations. First, for both algorithms, as described in Equations \eqref{eq:mp_def} and \eqref{eq:mp_optim_disc_def}, a \emph{full-search} has to be performed in $\mathcal{D}_{\rm d}=\dict(\Lambda_\ud)$ to compute all the squared scalar products $S_u$ of the current residue $u=R^m f$, with atoms $g_k$. We decide thus to reduce the computational complexity of this full-search with the help of the Fast Fourier Transform (FFT). One component (for 1-D signals) or two components (for 2-D signals) of the parametrization correspond indeed to a regular grid of atoms positions, which makes $S_u$ a discrete correlation relatively to these parameters.  
Moreover, as described in detail in \cite{Ventura2004,rosaphd}, we apply the fast \emph{boundary renormalization} of atoms, where atoms of $\mathcal{D}$ truncated by the limit of the signal remain valid atoms, i.e. of unit norm, and features that suddenly terminate at the signal boundary are correctly caught in the procedure. Notice that all our dMP and gMP experiments are performed in the non-weak case, i.e. $\alpha=1$. 

Second, for the Gradient-Ascent optimization, we realize some simplifications to the initial formulation~: the best discrete atom only is optimized at each MP iteration and $\kappa_m=\kappa>0$ for all $m$, with $\kappa$ typically equal to 5 or 10. Even if these two restrictions are not optimal compared to the method described in the theoretical results, the gain of the optimization in the quality of signals reconstructions is already impressive. We also set all the step sizes to $t_r=\chi>0$, with $\chi=0.1$ in all our experiments. Then, at each optimization step $r$, we adaptively decrease the step parameter $t_r$ by dividing it by 2 if the ascent condition is not met, i.e. if $S_u(\phi_{r+1}(k))< S_u(\phi_{r}(k))$. If after 10 divisions, the ascent condition still does not hold, the optimization process is simply stopped.

Finally, let us mention that our algorithms are written in MATLAB\copyright\ and are consequently not truly optimized. The different computation times that we provide through this section allow us only to compare various schemes, as for dMP and gMP decomposition of the same signal. All of our experiments were realized on a Pentium 1.73 GHz with 1Gb of memory.

\subsection{One Dimensional Analysis}
\label{subsec:1-D-experiments}

This section analyzes the benefit obtained from gMP, and from an increased density of the discrete dictionary, when decomposing some specific classes of randomly generated 1-D signals. In our experiments, each 1-D signal is of unit norm and has $N=2^{13}$ samples. Each signal consists of the sum of 100 random bursts, each burst being a rectangular or Gaussian window, depending on the class of the signal. The position and magnitude of each burst is selected at random, according to a uniform distribution. The duration of the rectangular window and the standard deviation of the Gaussian function are selected uniformly within the range $[\tinv{2}L, \frac{3}{4}L]$, for $L=2^8$. The mother function of the dictionary is the \emph{Mexican Hat} function $g(t)\propto(1-t^2)\,e^{-t^2/2}$. Its scale and translation parameters are sampled as defined in Section~\ref{subsec:example}, following the $\tau$-adic discretization $\Lambda_\ud = \{ (nb_0\tau^j,a_0\tau^j): j,n\in\Zbb\}$, with $a_0=1$. We work in the non-weak case, i.e. $\alpha=1$, for dMP and gMP, and we set $\kappa=10$ for gMP.

\begin{figure*}
  \centering
  \subfigure[Gaussian\label{fig:convergence-gauss}]{
    \psfrag{tag1}[l][][0.6]{\!\!\!\!\!dMP$(2,0.5)$}
    \psfrag{tag2}[l][][0.6]{\!\!\!\!\!gMP$(2,0.5)$}
    \psfrag{tag3}[l][][0.6]{\!\!\!\!\!dMP$(1,0.25)$}
    \psfrag{tag4}[l][][0.6]{\!\!\!\!\!gMP$(1,0.25)$}
    \psfrag{tagx}[][b][.8]{$m$}
    \psfrag{tagy}[][][.8]{$\|R^m f\|^2$}	
    \includegraphics[height=4cm,width=5.8cm]{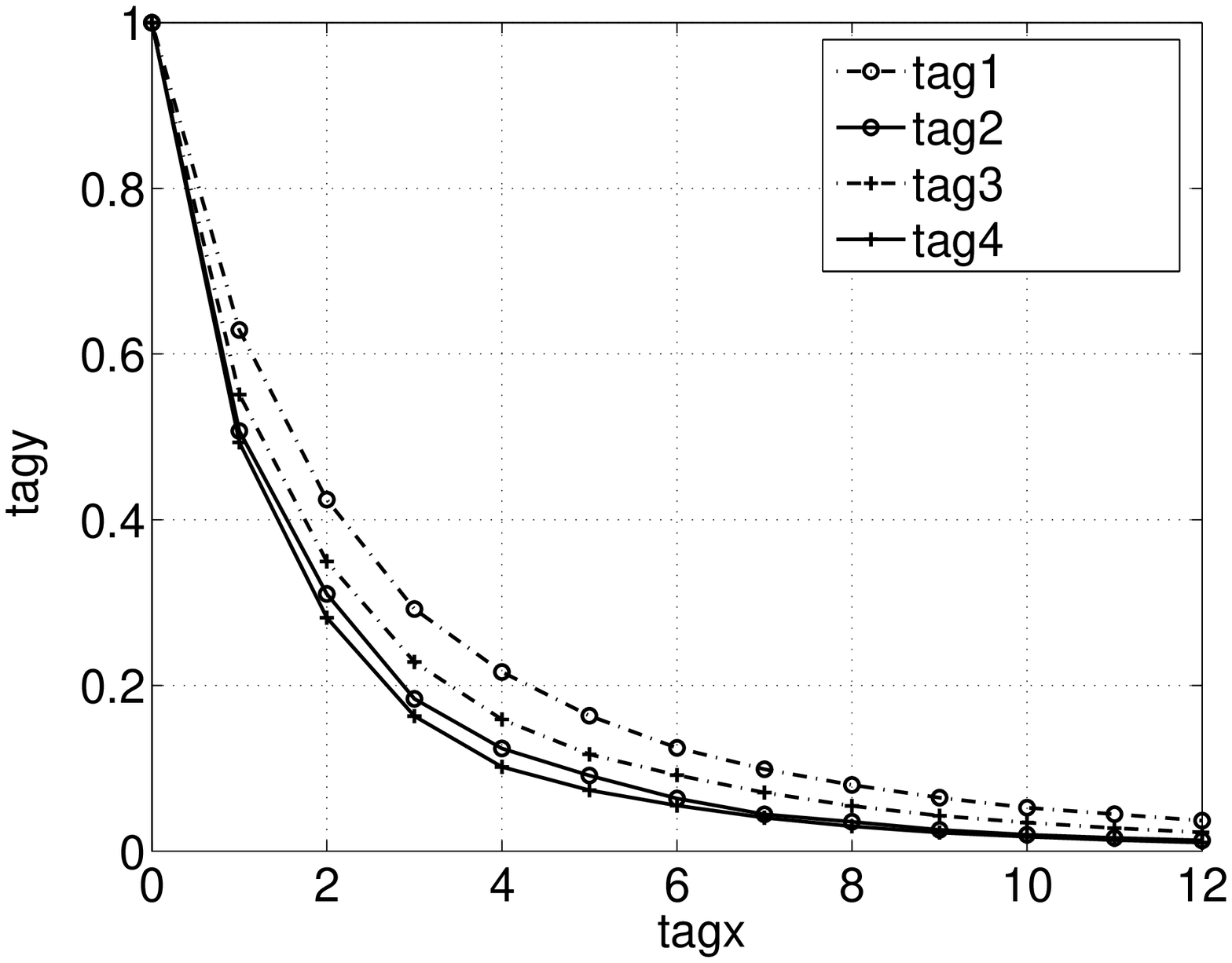}
  }
  \subfigure[Rectangular\label{fig:convergence-rect}]{
    \psfrag{tag1}[l][][0.6]{\!\!\!\!\!dMP$(2,0.5)$}
    \psfrag{tag2}[l][][0.6]{\!\!\!\!\!gMP$(2,0.5)$}
    \psfrag{tag3}[l][][0.6]{\!\!\!\!\!dMP$(1,0.25)$}
    \psfrag{tag4}[l][][0.6]{\!\!\!\!\!gMP$(1,0.25)$}
    \psfrag{tagx}[][b][.8]{$m$}
    \psfrag{tagy}[][][.8]{$\|R^m f\|^2$}	
    \includegraphics[height=4cm, width=5.8cm]{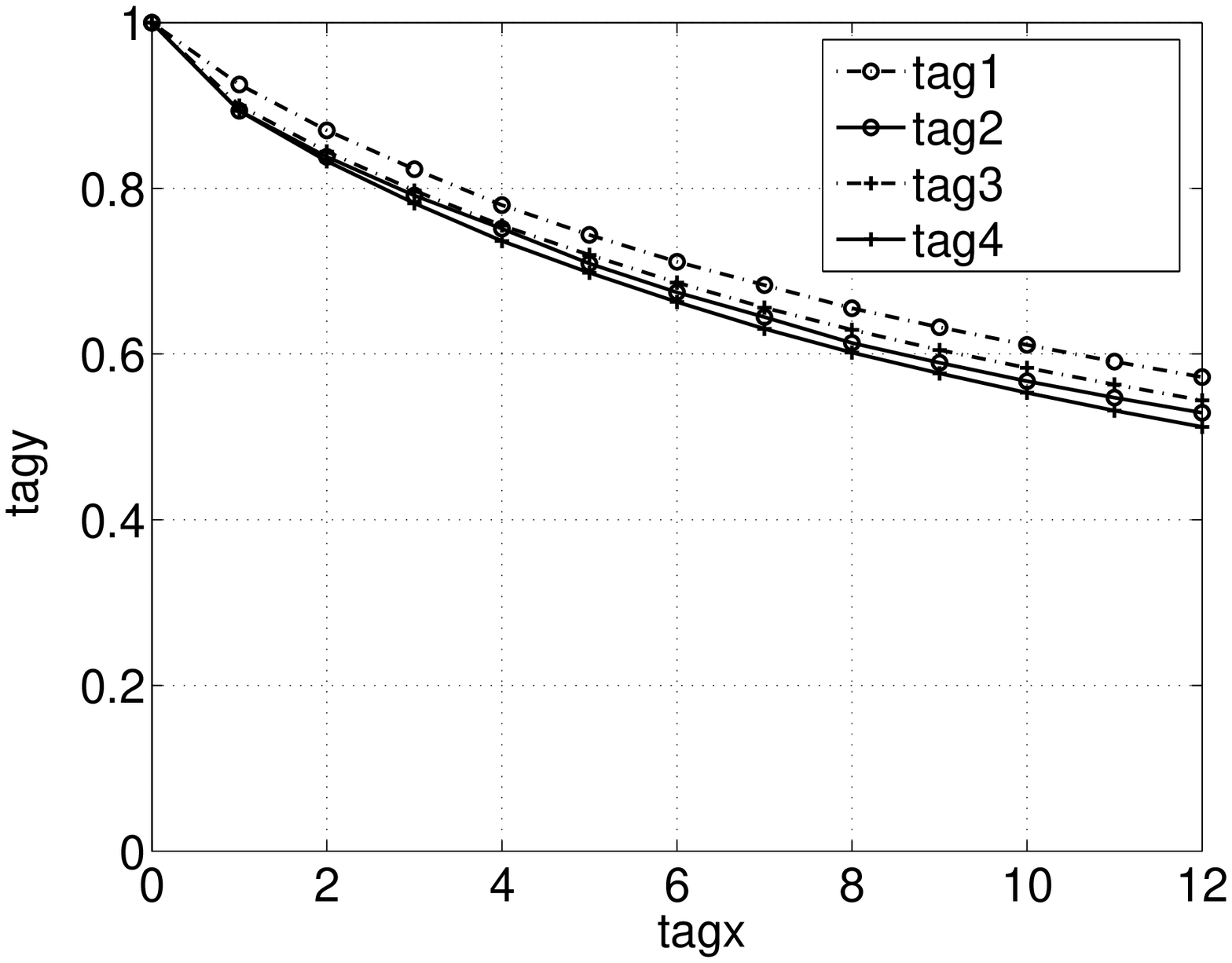}
  }\\
  \subfigure[Gaussian, $1^{st}$ MP step\label{fig:NAE-gauss-1}]{
    \psfrag{NAE}[][][.8]{NAE}
    \psfrag{dMP}[l][][.6]{\!\!\!\!\!dMP}
    \psfrag{oMP}[l][][.6]{\!\!\!\!\!gMP}
    \psfrag{a}[][][0.4]{$a$}
    \psfrag{tagx}[][][.9]{$\log_2\tau$}
    \includegraphics[keepaspectratio,height=4cm]{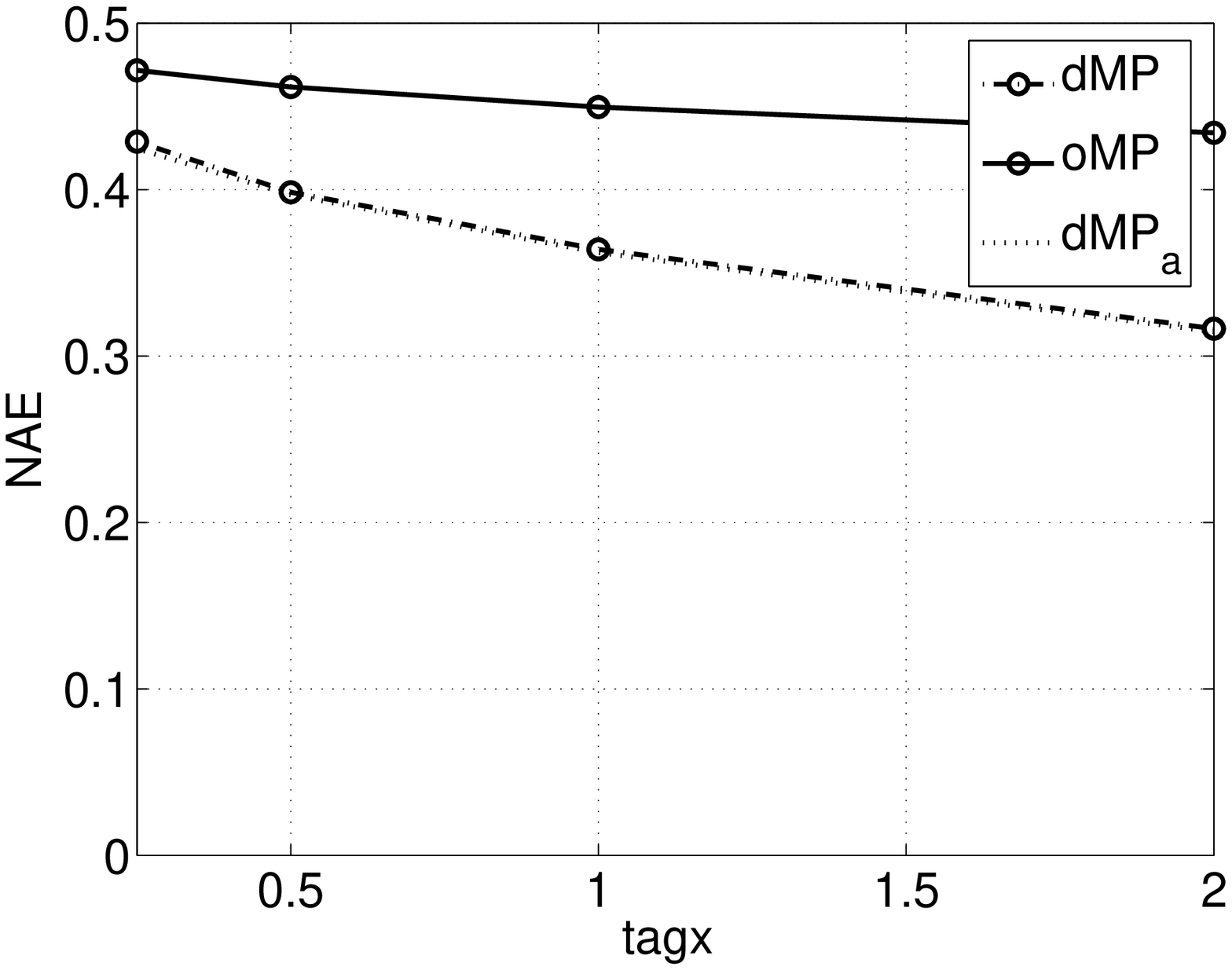}
  }
  \subfigure[Gaussian, $30^{th}$ MP step\label{fig:NAE-gauss-30}]{
    \psfrag{NAE}[][][.8]{NAE}
    \psfrag{dMP}[l][][.6]{\!\!\!\!\!dMP}
    \psfrag{oMP}[l][][.6]{\!\!\!\!\!gMP}
    \psfrag{a}[][][.4]{$a$}
    \psfrag{tagx}[][][.9]{$\log_2\tau$}
    \includegraphics[keepaspectratio,height=4cm]{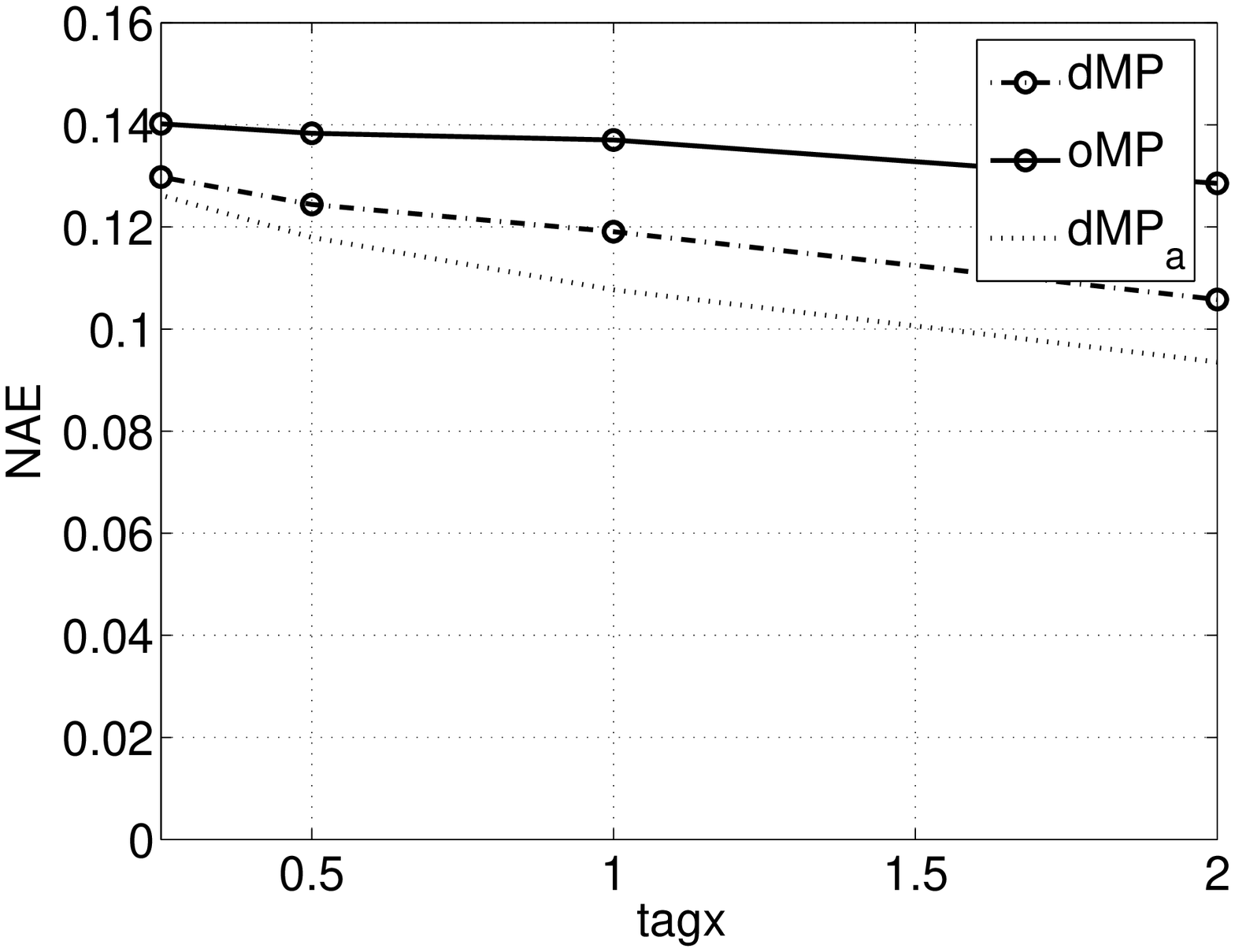}
  }\\
  \subfigure[Rectangular, $1^{st}$ MP step\label{fig:NAE-rect-1}]{
    \psfrag{NAE}[][][.8]{NAE}
    \psfrag{dMP}[l][][.6]{\!\!\!\!\!dMP}
    \psfrag{oMP}[l][][.6]{\!\!\!\!\!gMP}
    \psfrag{a}[][][.4]{$a$}
    \psfrag{tagx}[][][.9]{$\log_2\tau$}
    \includegraphics[keepaspectratio,height=4cm]{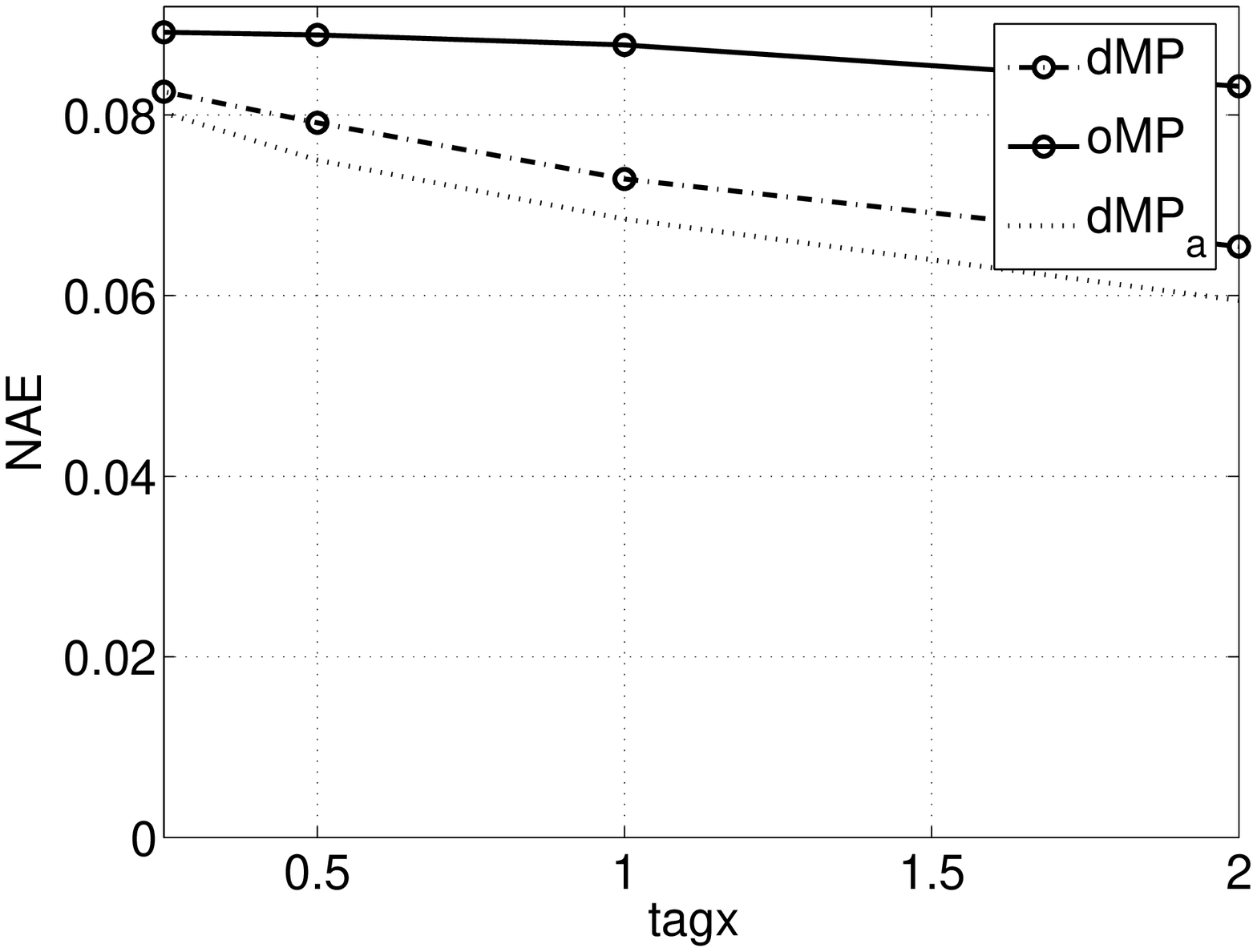}
  }
  \subfigure[Rectangular, $30^{th}$ MP step\label{fig:NAE-rect-30}]{
    \psfrag{NAE}[][][.8]{NAE}
    \psfrag{dMP}[l][][.6]{\!\!\!\!\!dMP}
    \psfrag{oMP}[l][][.6]{\!\!\!\!\!gMP}
    \psfrag{a}[][][.4]{$a$}
    \psfrag{tagx}[][][.9]{$\log_2\tau$}
    \includegraphics[keepaspectratio,height=4cm]{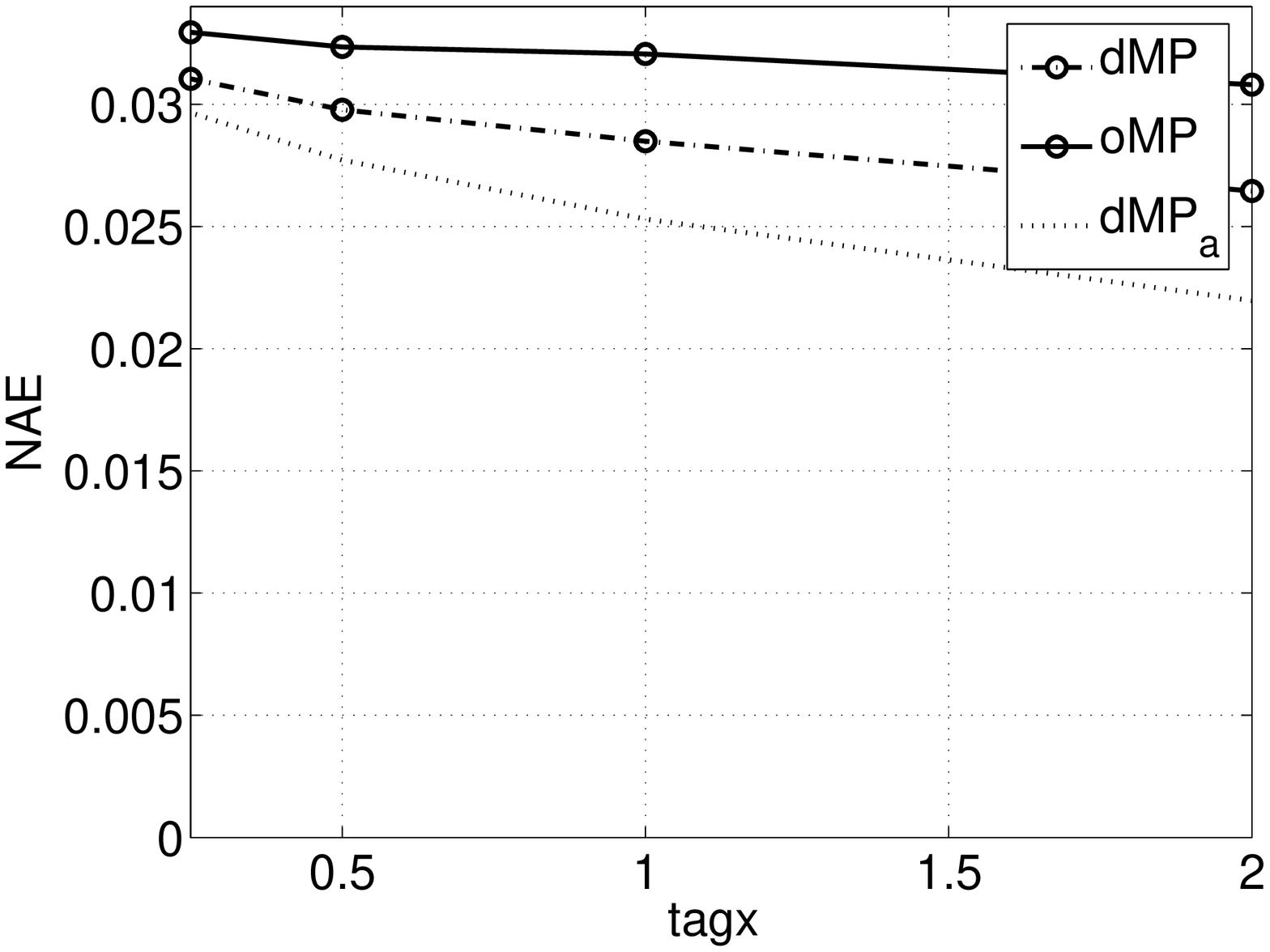}
  }
  \caption{(a)-(b) Residual energy as a function of the MP iteration. dMP$\,(b_0,\log_2\tau)$ and gMP$\,(b_0,\log_2\tau)$ refer to discrete and optimized MP, computed on a discretization $\Lambda_\ud = \{ (nb_0\tau^j,a_0\tau^j): j,n\in\Zbb\}$ of the continuous Mexican Hat dictionary. (c)-(f) Normalized atom energy (NAE) as a function of the $\log_2\tau$ discretization parameter. $b_0$ is set to one in all cases. dMP and gMP respectively refer to discrete and optimized MP. dMP$_a$ provides a lower bound to the decrease of NAE with $\log_2\tau$, and is formally described in the text.}
  \label{fig:convergence-NAE}
\end{figure*}

Figures~\ref{fig:convergence-gauss} and \ref{fig:convergence-rect} analyze how the energy $\|R^m f\|^2$ of the residual decreases with the number $m$ of MP iterations for the random Gaussian and rectangular signals, respectively. Notice that only a small number of iterations are studied (twelve) since our analysis aims at analyzing the behaviour of dMP and gMP on one class of signals. However the current residual $R^m f$ belongs only approximately to the considered class on small $m$ when not many atoms have been substracted to $f=R^0 f$. Results presented are averaged over 20 trials. In each graph, two distinct discretizations of the Mexican Hat parameters are considered to provide two discrete dictionaries, with one $(b_0=1, \log_2\tau = 0.25)$ being two times denser than the other $(b_0=2, \log_2\tau = 0.5)$, according to the behavior\footnote{Obviously equivalent for $\log_2\tau$ or $\ln\tau$ variations.} of the density radius $\rho_\ud$ analyzed in Section~\ref{subsec:example}. Both discrete and geometrically optimized MP are studied for each dictionary. We observe that gMP significantly outperforms dMP, and that an increased density of the dictionary also speeds up the MP convergence. By comparing Figure~\ref{fig:convergence-gauss} and \ref{fig:convergence-rect}, we also observe that the residual energy decreases much faster for Gaussian signals than for rectangular ones, which unsurprisingly reveals that the Mexican Hat dictionary is better suited to represent Gaussian structures. 

Figures~\ref{fig:NAE-gauss-1}-\ref{fig:NAE-rect-30} further analyze the impact of the discretization of the dictionary parameters on MP convergence. In these figures, we introduce the notion of {\it normalized atom energy} (NAE) to measure the convergence rate of a particular dictionary dealing with a specific class of signals at a specific MP iteration step. Formally, the NAE denotes the expected value of the best squared atom coefficient computed on a normalized signal when this one is randomly generated within a specific class of signals. Mathematically, ${\rm NAE} = \mathbb{E}\big[\scp{g_{\lambda_*}}{\tfrac{u}{\|u\|}}^2\big]$, where $u$ is a sample signal of the class and the $g_{\lambda_*}$ the associated best atom for a fixed greedy algorithm (dMP or gMP). We show the dependence of NAE on the discretization for the $1^{st}$ and $30^{th}$ iteration\footnote{Note that the NAE at the $30^{th}$ iteration refers to the NAE computed on the residual signals obtained after 29 iterations of the gMP with the densest dictionary, independently of the actual discrete dictionary considered at iteration 30. Hence, the reference class of signals to compute the NAE at iteration 30 is the same for all investigated dictionaries, i.e. for all $\log_2\tau$ values.} for both rectangular and Gaussian signals. Results are averaged over 500 trials. 

By considering the dMP and gMP curves in Figures~\ref{fig:NAE-gauss-1}-\ref{fig:NAE-rect-30}, we first observe that the NAE is significantly higher for gMP than for dMP, which confirms the advantage of using gradient ascent optimization to refine the parameters of the atoms extracted by dMP. Note that the NAE for a Gaussian random signal (Fig. \ref{fig:NAE-gauss-1}-\ref{fig:NAE-gauss-30}) is nearly one order of magnitude higher than for a rectangular one (Fig.\ref{fig:NAE-rect-1}-\ref{fig:NAE-rect-30}). This confirms that the Mexican Hat dictionary better matches the Gaussian structures than the rectangular ones. In all cases, the NAE sharply decreases with the iteration index, which is not a surprise as the coherence between the signal and the dictionary decreases as MP expansion progresses.

To better understand the penalty induced by the discretization of the continuous dictionary, we now analyze how the rate of convergence for a particular class of signals behaves compared to the reference provided by a signal composed of a single Mexican Hat function. For that purpose, an additional curve, denoted dMP$_a$, has been plotted in each graph. This curve is expected to provide an upper bound to the penalty induced by a sparser dictionary. Specifically, dMP$_a$ plots the energy captured during the $1^{st}$ step of the dMP expansion of a random (scale and position) Mexican Hat function, as a function of the discretization parameter $\log_2\tau$. As the Mexican Hat is the generative function of the dictionary, the $1^{st}$ step of the MP expansion would capture the entire function energy if the entire continuous dictionary were used, but is particularly penalized by a discretization of the dictionary. In each graph of Figures~\ref{fig:NAE-gauss-1}-\ref{fig:NAE-rect-30}, to compare dMP$_a$ with dMP, the dMP$_a$ curve obtained with pure atoms (i.e. unit coefficients) is scaled to correspond to atoms whose energy is set to the NAE expected from the expansion of the corresponding class of signals with a continuous dictionary. In practice, the NAE expected with a continuous dictionary is estimated based on the NAE computed with gMP and the densest dictionary ($\log_2\tau=0.25$). The approximation is reasonable as we observe that gMP saturates for small $\log_2\tau$ values, i.e. for large densities. We first observe that both the dMP and the dMP$_a$ curves nearly coincide in Figure~\ref{fig:NAE-gauss-1}. Hence, the MP expansion of a Gaussian signal is penalized as much as the one of a Mexican Hat function by a reduction of the discrete dictionary density. We then observe that the penalty induced by a reduction of density decreases as the coherence between signal and dictionary structures drops. This is for example the case when the signal to represent is intrinsically sharper than the dictionary structures (Fig.~\ref{fig:NAE-rect-1}-\ref{fig:NAE-rect-30}), or because the coherent structures have been extracted during the initial MP steps (Fig.~\ref{fig:NAE-gauss-30}). This last observation is of practical importance because it reveals that using a coarsely discretized dictionary incurs a greater penalty during the first few iterations of the MP expansion than during the subsequent ones. For compression applications, it might thus be advantageous to progressively decrease the density of the dictionary along the expansion process, the cost associated to the definition of the atom indices decreasing with the density of the dictionary\footnote{Less distinct atom indices need to be described by the codewords.}. Hence, it might be more efficient -- in a rate-distortion sense -- to use a dense but expensive dictionary during the first MP iterations, so as to avoid penalizing the MP convergence rate, but a sparser and cheaper during subsequent steps, so as to save bits. We plan to investigate this question in details in a future publication.

\subsection{Two Dimensional Analysis}
\label{subsec:2-D-experiments}

This section analyzes experimentally the effect of discretizing a dictionary on the Matching Pursuit decomposition of images, i.e. with the Hilbert space $L^2(\Rbb^2)$. 

\paragraph{Parametrization and Dictionary } We use the same dictionary as in \cite{Ventura2004}. Its mother function $g$ is defined by a separable product of two 1-D behaviors : a Mexican Hat wavelet in the $x$-direction, and a Gaussian in the $y$-direction, i.e. $g(\vec{x}) = (\tfrac{4}{3\pi})^{1/2}\,(1 - x^2)\,\exp(-\tinv{2}\,|\vec{x}|^2)$, where $\vec{x}=(x,y)\in \Rbb^2$ and $\|g\|=1$ \cite{mall98}. Notice that $g$ is infinitely differentiable. 

The dictionary is defined by the translations, rotations, and
anisotropic dilations of $g$. Mathematically, these transformations are
represented by operators $T_{\vec{b}}$, $R_\theta$, and
$D_{\vec{a}}$, respectively. These are given by $[T_{\vec{b}}\,g]\big(\vec{x}\big) = g\big(\vec{x}-\vec{b}\big)$, $[R_\theta\,g]\big(\vec{x}\big) = g\big(r^{-1}_\theta\,\vec{x}\big)$, and $[D_{\vec{a}}\,g]\big(\vec{x}\big) = (a_1a_2)^{-1/2}\,g\big(d^{\,-1}_{\vec{a}}\vec{x} \big)$, for $\theta\in S^1 \simeq [0,2\pi)$, $\vec{b}\in \Rbb^2$, $\vec{a}=(a_1,
a_2)$, $a_1, a_2 \in \Rbb^*_+$, while $r_\theta$ is the usual $2\times2$ rotation matrix $r_\theta$ and $d_{\vec{a}} = {\rm diag}(a_1,a_2)$.

In other words, we have a parametrization of $P=5$ dimensions and
$\Lambda=\{\lambda=(\lambda^0,\ldots,\lambda^4)=(b_1,b_2,\theta,a_1,a_2)\in
\Rbb^2\times S^1\times(\Rbb^*_+)^2\}$.
At the end, each atom of the dictionary
$\mathcal{D}=\{g_\lambda:\lambda\in\Lambda\}$ is generated by $g_\lambda(\vec{x}) = [U(\lambda)\,g](\vec{x}) \triangleq [T_{\vec{b}}\,R_{\theta}\,D_{\vec{a}}\,g](\vec{x})$, with $\|g_\lambda\|=\|g\|=1$.

Obviously, the dictionary $\mathcal{D}$ is complete in $L^2(\Rbb^2)$. Indeed, translations,
rotations and isotropic dilations alone are already enough to constitute a wavelet
basis of $L^2(X)$ since $g$ is an admissible wavelet
\cite{daub92,jpamurvdgali03}. Finally, as requested in the previous section,
from the smoothness of $g$ and of the transformations $U$ above, the atoms
$g_\lambda$ of our dictionary $\mathcal{D}$ are twice
differentiable on each component $\lambda^i$.

\paragraph{Spatial Sampling} For all our experiments, images are discretized on a Cartesian regular
grid of pixels, i.e. an image $f$ takes its values on the grid
$\mathcal{X}=\big([0,N_x)\times[0,N_y)\big)\cap\Zbb^2$, with
$N_x$ and $N_y$ the ``$x$'' and ``$y$'' sizes of the grid. We work in the \emph{continuous approximation}, that is we assume that the grid $\mathcal{X}$ is fine enough to guarantee that the scalar products $\scp{\cdot}{\cdot}$ and norms $\|\cdot\|$ are well estimated from their discrete counterparts. 
This holds of course for band-limited functions on $L(\Rbb^2)$.

In consequence, in order to respect this continuous approximations and to have dictionary atoms smaller than the image size, the mother function $g$ of our dictionary $\mathcal{D}$ must be dilated in a particular range of scales so that $g_\lambda$ is essentially band-limited, i.e. $a_1,a_2\in [a_{\rm m},a_{\rm M}]$. According to the definition of $g$ above, we set experimentally $a_{\rm m}=0.7$ and $a_{\rm M}=\min(N_x,N_y)$.

\paragraph{Discrete Parameter Space}

We decide to sample regularly $\Lambda$ so that to have $N_{\rm pix}=N_xN_y$ positions $\vec{b}$, $J^2$ scales $a_1$ and $a_2$ selected logarithmically in the range $[a_{\rm m},a_{\rm M}]$, and $K$ orientations evenly spaced in $[0, \pi)$ , with $J,K\in\Nbb$. At the end, we obtain the discretized parameter set $\Lambda_{\rm d} = \Lambda_{\rm d}(N_{\rm pix},J,K) =\big\{\,(\vec{b},\theta_n,a_{1j},a_{2j'}),\ \vec{b}\in\mathcal{X},\ n\in[0,K-1],\ j,j'\in[0,J-1]\,\big\},$ and the corresponding dictionary $\mathcal{D}_{\rm d}(N_{\rm
  pix},J,K)=\dict(\Lambda_{\rm d}(N_{\rm pix},J,K))$.  The number of atoms in the dictionary is simply
$|\mathcal{D}_{\rm d}|=J^2K\,N_{\rm pix}$.

\begin{figure}
\centering
  \subfigure[\label{fig:barbara-no-MP-J5K8}]{
    \includegraphics[keepaspectratio,width=5cm]{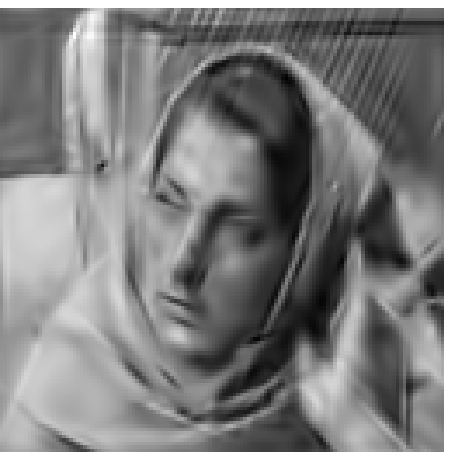}
  }
  \subfigure[\label{fig:barbara-op-MP-r10}]{
    \includegraphics[keepaspectratio,width=5cm]{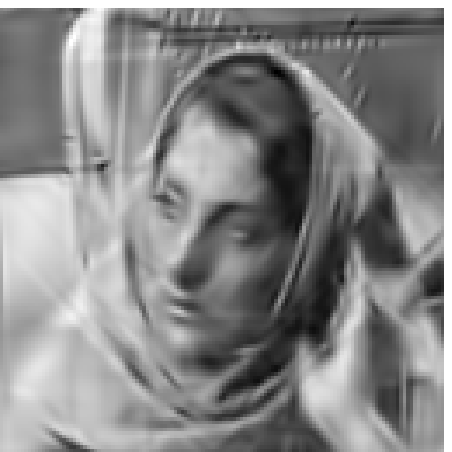}
  }
  \caption{300 atoms reconstruction results. (a) dMP : $J=5$, $K=8$, PSNR\,: $26.63$\,dB, 4634s. 
(b) gMP : $J=3$, $K=4$, $\kappa=10$, PSNR\,: $26.68$\,dB, 949s.}
\end{figure}

\begin{table}
\centering
\begin{tabular}{|c|l|l|}
\hline
&\multicolumn{1}{|c|}{$J = 3$}&\multicolumn{1}{|c|}{$J = 5$}\\
\hline
&&\\[-2.5mm]
$K = 4$&24.30\,dB \phantom{1}(834s)&25.88\,dB (2327s)\\
\cline{2-3}
&&\\[-2.5mm]
$(\kappa=5)$&26.08\,dB \phantom{1}(889s)&27.09\,dB (2381s)\\
$(\kappa=10)$&26.68\,dB \phantom{1}(950s)&27.37\,dB (2447s)\\
\hline
&&\\[-2.5mm]
$K = 8$&25.21\,dB (1660s)&26.63\,dB (4634s)\\
\cline{2-3}
&&\\[-2.5mm]
$(\kappa=5)$&27.05\,dB (1715s)&27.92\,dB (4703s)\\
$(\kappa=10)$&27.44\,dB (1772s)&28.09\,dB (5131s)\\
\hline
\end{tabular}
\caption{{\rm d}MP and {\rm g}MP applied on {\tt Barbara} image. Quality (in
PSNR) of the reconstruction after 300 iterations for various $J$, $K$ and $\kappa$. In each table cell, the first row correspond to {\rm d}MP result, the second and the third rows to {\rm g}MP.}\ \\[-5mm]
\label{tab:dmp-omp-compar}
\end{table}

\paragraph{Results }

We start our experiment by decomposing the venerable image of {\tt Barbara}. 300 atoms were selected by dMP and gMP for various $J$ and $K$. Results are presented in Table \ref{tab:dmp-omp-compar}. In these tests, the best quality obtained for dMP corresponds obviously to the finest grid, i.e. $J=5$ and $K=8$ (26.63\,dB, Fig.\ref{fig:barbara-no-MP-J5K8}), with a computational time (CT) of 4634s. With 10 optimization steps ($\kappa=10$), the gMP for the coarsest parametrization ($J=3$ and $K=4$) is equivalent to the best dMP result with a PSNR of 26.68\,dB and a CT of only 950s, i.e. almost five time faster. This is also far better than the dMP on the same grid (24.30\,dB). 
The visual inspection of the dMP image ($J=5$, $K=8$, Fig.\ref{fig:barbara-no-MP-J5K8}) and the gMP image ($J=3$, $K=4$, $\kappa=10$, \ref{fig:barbara-op-MP-r10}) is also instructive. Most of the features of the gMP results are well represented (e.g. Barbara's mouth, eyes, nose, hair, ...). However, the regular pattern of the chair in the background of the picture, which needs a lot of similar atoms, is poorly drawn. This can be explained by the fact that this highly directional structure has to be represented by a lot of similarly oriented and scaled atoms with similar amplitude. The fine grid of dMP has therefore more chance to correctly fit these atoms, while the gMP on its coarse grid is deviated in its optimization process to more prominent structure with higher amplitudes. Notice finally, the best optimized result (PSNR $28.09$\,dB) is obtained for $\kappa=10$ on the grid associated to $J=5$ and $K=8$ orientations.

\medskip
For our second experiment, we compare dMP and gMP ($\kappa=10$) 300 atoms approximation of well known 128$\times$128 pixels pictures, namely {\tt Lena}, {\tt Baboon}, {\tt Cameraman}, {\tt GoldHill}, and {\tt Peppers}, on the same pa\-ra\-me\-tri\-za\-tion grid ($J=4$, $K=8$). For a computational time slightly higher (5\%) than the dMP decomposition, we reach in all cases a significantly higher PSNR with gMP than with dMP, i.e. the dB gain is within the range $[0.87, 2.03]$.

\begin{table}
  \centering
  \begin{tabular}{|r|c|c|}
    \hline
    Image name&dMP&gMP ($\kappa=10$)\\
    \hline
&&\\[-2.5mm]
    {\tt Barbara}&25.94\,dB (2707s)&27.86\,dB (2820s)\\
    {\tt Lena}&26.50\,dB (2709s)&28.53\,dB (2857s)\\
    {\tt Baboon}&24.06\,dB (2770s)&24.93\,dB (2900s)\\
    {\tt Cameraman}&25.80\,dB (2807s)&27.62\,dB (2918s)\\
    {\tt GoldHill}&26.54\,dB (2810s)&28.12\,dB (2961s)\\
    {\tt Peppers}&24.51\,dB (2853s)&26.69\,dB (3013s)\\
    \hline
  \end{tabular}
  \caption{Comparison of {\rm d}MP and {\rm g}MP on
    different usual images of size $128\!\times\!128$. Computations have been performed for $J=4$, $K=8$, 300 atoms. Computation times are given indicatively in parenthesis.
}\ \\[-10mm]
  \label{tab:diff-images-FSMP-nonopt-vs-opt}
\end{table}

\section{Related Works}
\label{sec:rel-work}
A similar approach to our geometric analysis of MP atom selection rule has been proposed in \cite{watson1998sai}. In that paper, a dictionary of ($L^2$-normalized) wavelets is seen as a manifold associate to a Riemannian metric. However, the authors restrict their work to wavelet parametrization inherited from Lie group (such as the affine group). They also work only on the $L^2$ (dictionary) distance between dictionary atoms and do not introduce intrinsic geodesic distance. They define a discretization of the parametrization $\Lambda$ such that, in our notations, $\mathcal{G}_{ij}\Delta\lambda^i\Delta\lambda^j < \epsilon$, with $\Delta\lambda(k)$ the local width of the cell localized on $k\in\Lambda_{\rm d}$. There is however no analysis of the effect of this discretization on the MP rate of convergence.

In \cite{gribonval2001fmp}, the author uses a 4-dimensional Gaussian chirp dictionary to analyze 1-D signals with MP algorithm. He develops a fast procedure to find the best atom of this dictionary in the representation of the current MP residual by applying a two-step search. First, by setting the chirp rate parameter to zero, the best common Gabor atom is found with full search procedure taking advantage of the FFT algorithm. Next, a ridge theorem proves that starting from this Gabor atom, the best Gaussian chirp atom can be approximated with a controlled error. The whole method is similar to the development of our optimized matching pursuit since we start also from a discrete parametrization to find a better atom in the continuous one. However, our approach is more general since we are not restricted to a specific dictionary. We use the intrinsic geometry of any smooth dictionary manifold to perform a optimization driven by a geometric gradient ascent. 

\section{Conclusions}
\label{sec:conclusion}

In this paper, we have adopted a geometrical framework to study the effect of dictionary discretization on the rate of convergence associated to MP. In a first step, we have derived an upper bound for this rate using geometrical quantities inherited from the dictionary seen as a manifold, such as the geodesic distance, the condition number of the dictionary, and the covering property of the discrete set of atoms in the continuous dictionary. We have also shown in a second step how a simple optimization of the parameters selected by the discrete dictionary, can lead theoretically and experimentally to important gain in the approximation of (general) signals. 

In a future study, it could be interesting to see how our methods extend to other greedy algorithms, like the Orthogonal Matching Pursuit (OMP) \cite{prk93}. However, this extension has to be performed carefully since we need to characterized the convergence of continuous OMP, as it is here for the one of MP induced by the existence of a greedy factor. 

Our work paves the way for future extensions and advances in two practical fields. As explained in our 1-D experiments, a first idea could be to analyze carefully the benefit -- in a rate-distortion sense -- of using a dense but expensive dictionary during the first MP iterations, so as to avoid penalizing the MP convergence rate, but a sparser and cheaper dictionary during subsequent steps, so as to save bits. We plan to investigate this question in details in a future publication.

Another idea is to analyze the behaviors of gMP in the Compressive Sensing (CS) formalism, that is after random projection of the signal and atoms. Matching Pursuit is already used currently as a retrieval algorithm of CS of sparse signals \cite{candes2006rup,duarte2006ssd,CSredunddico}. However, recent results \cite{randprojmanif} suggests also that for manifold of bounded condition number, their geometrical structure (metric, distances) is essentially preserved after random projection of their points in a smaller space than the ambient one. If a natural definition of random projection in our continuous formalism can be formulated, a natural question is thus to check if the gradient ascent technique survives after random projection of the residual and the atoms on the same subspace. This could lead to dramatic computation time reduction, up to controlled errors that could be even attenuated by the greedy iterative procedure.

\section*{Acknowledgements}
LJ wishes to thank Prof. Richard Baraniuk and his team at Rice University (Houston, TX, USA) for the helpful discussions about general ``manifolds processing'' and Compressive Sensing. LJ is also very grateful to R. Baraniuk for having accepted and funded him during a short postdoctoral stay at Rice University. LJ and CDV are funded by the Belgian FRS-FNRS. We would like to thank Dr. David Kenric Hammond (LTS2/EPFL, Switzerland) for his careful proofreading and the referees for valuable comments on this paper.

\appendix 

\section{Complements on the Geometry of $(\Lambda,\mathcal{G}_{ij})$}
\label{app:comp-diff-geo}

In this short appendix, we provide some additional information on the geometrical concepts developed in Section \ref{sec:dico-param-manifold}. First, as explained in that section, the parameter space $\Lambda$ of the dictionary $\mathcal{D}=\dict(\Lambda)$ is linked to a Riemannian manifold $\mathcal{M}=(\Lambda,\mathcal{G}_{ij})$ with a structure inherited from the dictionary $\mathcal{D}\subset L^2(X)$. From the geodesic definition \eqref{eq:geodesic-distance-dico-def} and the metric relation \eqref{eq:pullback-relation-fundation}, we see that the curve $\gamma_{_{\lambda_a\lambda_b}}(t)\in\Lambda$ is thus also a geodesic in $\mathcal{M}$. In other words, it is defined only from the metric $\mathcal{G}_{ij}$ and not anymore from the full behavior of atoms of $\mathcal{D}\subset L^2(X)$. In \cite{TRLJ0701}, we explain also that $\mathcal{M}$ is in fact an \emph{immersed manifold} \cite{carmo1992rg} in the Hilbert manifold $\mathcal{D}\subset L^2(X)$, and $\mathcal{G}_{ij}$ is the associated \emph{pullback} metric. All the geometric quantities of the Riemannian analysis of $\mathcal{M}$, such as Christoffel's symbols, covariant derivatives, curvature tensors, etc. can be defined. This is actually done in the following appendices of this paper. 

Second, some important designations can be introduced. The metric $\mathcal{G}_{ij}(\lambda)$ is a (\emph{covariant}) \emph{tensor} of rank-2, i.e. described by two subscript indices, on $\mathcal{M}$. This means that $\mathcal{G}_{ij}$ satisfies a specific transformation under changes of coordinates in $T_\lambda\Lambda$ such that the values of the bilinear form\footnote{Also named first fundamental form \cite{carmo1992rg}.} $\mathcal{G}_\lambda(\xi,\zeta) \triangleq \xi^i\,\zeta^j\,\mathcal{G}_{ij}(\lambda)$ that it induces are unmodified\footnote{In the same way that the scalar product between two vectors in the usual Euclidean space is independent of the choice of coordinates.}. A function $f:\Lambda\to \Rbb$ is a \emph{scalar field} on $\mathcal{M}$, or rank-0 tensor. A vector field $\zeta^i(\lambda)$ on this manifold, which associates to each point $\lambda$ a vector in the tangent plane $T_{\lambda}\Lambda$, is a function $\zeta:\Lambda\to T_{\lambda}\Lambda\simeq\Rbb^P$ also named (\emph{contravariant}) rank-1 tensor, i.e. with one superscript. More generally, a rank-$(m,n)$ tensor is a quantity $T^{i_1\,\cdots\,i_m}_{j_1\,\cdots\,j_n}(\lambda)$ $m$-times contravariant and $n$-times covariant 
such that $\mathcal{G}_{i_1k_1}\cdots\,\mathcal{G}_{i_mk_m}\ \xi_1^{k_1}\cdots\,\xi_m^{k_m}\ T^{i_1\,\cdots\,i_m}_{j_1\,\cdots\,j_n}(\lambda)\ \zeta_1^{j_1}\cdots\,\zeta_n^{j_n}$ is invariant under change of coordinates in $T_\lambda\Lambda$ for any vectors $\{\xi_1,\,\cdots,\xi_m,\zeta_1,\,\cdots,\zeta_n\}$ in this space.

\section{Proof of Proposition \ref{prop:bounds-dico-curv}}
\label{app:proof-prop-unit-lower-bound-dico-curv}

Let $\gamma$ be a geodesic in $\mathcal{M}$ with curvilinear parametrization, i.e. with $|\gamma'(s)|=1$. Writing $\gamma=\gamma(s)$ and $\gamma'=\tfrac{\ud}{\ud s}\gamma(s)$, we have $\tfrac{\ud}{\ud s}\,g_{\gamma(s)} = \partial_i g_{\gamma}\,{\gamma'}^i$ and $\tfrac{\ud^2}{\ud s^2}\,g_{\gamma(s)} = \partial_{ij} g_{\gamma}\,{\gamma'}^i{\gamma'}^j + \partial_k g_{\gamma}\,{\gamma''}^k$, where we write abusively $\partial_i g_{\gamma}=\partial_i g_{\lambda}|_{\lambda=\gamma(s)}$ and similarly for second order derivative.

We need now some elements of differential geometry. Since $\gamma$ is a geodesic in $\mathcal{M}$, it respects the second order differential equation
$ {\gamma''}^k + \Gamma^k_{ij}\,{\gamma'}^i{\gamma'}^j = 0$, where the values $\Gamma^k_{ij} =  \tinv{2}\,\mathcal{G}^{lk}\,\big(\partial_j\,\mathcal{G}_{li} + \partial_i\,\mathcal{G}_{jl} - \partial_l\,\mathcal{G}_{ij}\big)$ are the Christoffel's symbols \cite{carmo1992rg} derived from the metric $\mathcal{G}_{ij}$. Therefore, we get 
\begin{eqnarray}{rcl}
\label{eq:second-s-der-gl-christ}
\tfrac{\ud^2}{\ud s^2}\,g_{\gamma}&=&\partial_{ij} g_{\gamma}\,{\gamma'}^i{\gamma'}^j - \partial_k g_{\gamma}\,\Gamma^k_{ij}\,{\gamma'}^i{\gamma'}^j\\
\label{eq:second-s-der-gl-cov-diff}
&=&\nabla_{ij} g_{\gamma}\,{\gamma'}^i{\gamma'}^j,
\end{eqnarray}
where $\nabla_i g_{\gamma}=\partial_i g_{\gamma}$ and $\nabla_{ij} g_{\gamma}=\nabla_i\nabla_j g_{\gamma} = \partial_{ij} g_{\gamma} - \partial_k g_{\gamma}\,\Gamma^k_{ij}$ are by definition the first order $i$ and the second order $ij$ \emph{covariant derivatives} of $g_{\gamma}$ respectively \cite{carmo1992rg}. In addition, we can easily compute that for $\mathcal{M}=(\Lambda,\mathcal{G}_{ij})$, 
\begin{equation}
\label{eq:christoffel-simplif-D}
\Gamma^k_{ij}\ =\ \mathcal{G}^{kl}\,\scp{\partial_{ij}g_{\lambda}}{\partial_lg_\lambda}.
\end{equation}

The lower bound of the proposition comes simply from the projection of $\tfrac{\ud^2}{\ud s^2}\,g_{\gamma(s)}$ onto $g_{\gamma}$. Indeed, for any $\lambda\in\Lambda$, since $\|g_\lambda\|^2=\scp{g_\lambda}{g_\lambda}=1$, $\scp{\partial_i g_\lambda}{g_\lambda}=0$ and $\scp{\partial_{ij} g_\lambda}{g_\lambda}=-\mathcal{G}_{ij}$. By \eqref{eq:second-s-der-gl-christ}, $\scp{\tfrac{\ud^2}{\ud s^2}\,g_{\gamma(s)}}{g_{\gamma}} = \scp{\partial_{ij}g_{\gamma}}{g_{\gamma}}\,{\gamma'}^i{\gamma'}^j = -\mathcal{G}_{ij}\,{\gamma'}^i{\gamma'}^j = -1$, and using Cauchy-Schwarz we get $\|\tfrac{\ud^2}{\ud s^2}\,g_{\gamma(s)}\| \geq 1$. Therefore, for $\epsilon>0$ and $\gamma_{\xi}~:~[0,\epsilon] \to \Lambda$, a segment of geodesic starting from $\lambda$ with unit speed $\xi$, $$\mathcal{K}\ \geq\ \sup_{\xi: |\xi|=1} \|\tfrac{\ud^2}{\ud s^2}\,g_{\gamma_\xi(s)}\big|_{s=0}\|\ \geq\ 1.$$

For the upper bound, coming back to any geodesic $\gamma$, we need to analyze directly $\|\tfrac{\ud^2}{\ud s^2}\,g_{\gamma(s)}\|^2$. Using \eqref{eq:second-s-der-gl-cov-diff} and the expression \eqref{eq:christoffel-simplif-D} of the Christoffel's symbols above, we have
$\|\tfrac{\ud^2}{\ud s^2}\,g_{\gamma(s)}\|^2$ $=$ $\| \nabla_{ij}\, g_{\gamma}\,{\gamma'}^i{\gamma'}^j\|^2$ $\leq \scp{\nabla_{ij}\, g_{\gamma}}{\nabla_{kl}\, g_{\gamma}}\,\mathcal{G}^{ik}\,\mathcal{G}^{jl}$, where we used $|{\gamma'}|=1$ and the Cauchy-Schwarz (CS) inequality expressed in the Einstein's summation notation on rank-2 tensors. This latter states that, for the tensors $A_{ij}=\nabla_{ij}\, g_{\gamma}$ and $B^{ij}={\gamma'}^i{\gamma'}^j$, $|A_{ij}B^{ij}|^2 \leq  |A_{ij}\,A_{kl}\,\mathcal{G}^{ki}\,\mathcal{G}^{lj}|\,|B^{ij}\,B^{kl}\,\mathcal{G}_{ki}\,\mathcal{G}_{lj}|$, the equality holding if the two tensors are multiple of each other. We prove in \cite{TRLJ0701} the general explanation for rank-n tensor as a simple consequence of the positive-definiteness of $\mathcal{G}_{ij}$. 

Therefore, taking $\gamma=\gamma_{\xi}$, and since $\gamma_\xi(0)=\lambda$,
\begin{equation}
\label{eq:covariant-KD-upper-bound}
\mathcal{K}\ \leq\ \sup_{\lambda\in\Lambda}\, \left[\bscp{\nabla_{ij}\,g_{\lambda}}{\nabla_{kl}\,g_{\lambda}}\,\mathcal{G}^{ik}\,\mathcal{G}^{jl}\,\right]^{\tinv{2}}.
\end{equation}

In the companion Technical Paper \cite{TRLJ0701}, we prove that this inequality is also equivalent to 
$$
\mathcal{K}\ \leq\ \sup_{\lambda\in\Lambda}\ \big[\,R(\lambda) + \|\Delta\,g_\lambda\|^2\,\big]^{\tinv{2}},
$$ 
where $R$ is the scalar curvature of $\mathcal{M}$, i.e. the quantity $R=R_{ijkl}\,\mathcal{G}^{ik}\,\mathcal{G}^{jl}$ contracted from the curvature tensor $R_{iklm} = \tinv{2}(
\partial_{kl}\mathcal{G}_{im} + \partial_{im} \mathcal{G}_{kl} - \partial_{km} \mathcal{G}_{il} - \partial_{il} \mathcal{G}_{km})
+ \mathcal{G}_{np} (\Gamma^n_{kl}\,\Gamma^p_{im} - \Gamma^n_{km}\,\Gamma^p_{il})$, and $\Delta g_\lambda=\mathcal{G}_{ij}\,\nabla^i\,\nabla^j\,g_\lambda$ is the Laplace-Beltrami operator applied on $g_\lambda$. The curvature $R$ requires only the knowledge of $\mathcal{G}_{ij}(\lambda)$ (and its derivatives), implying just one step of scalar products computations, i.e. integrations in $L^2(X)$.

The reader who does not want to deal with differential geometry can however get rid of the covariant derivatives of Equation \eqref{eq:covariant-KD-upper-bound} by replacing them by usual derivatives. This provides however a weaker bound. Indeed, using the expression \eqref{eq:christoffel-simplif-D} of the Christoffel's symbols, some easy calculation provides 
$0 \leq \bscp{\nabla_{ij}\,g_{\lambda}}{\nabla_{kl}\,g_{\lambda}}\,\mathcal{G}^{ik}\,\mathcal{G}^{jl} = \bscp{\partial_{ij}\,g_{\lambda}}{\partial_{kl}\,g_{\lambda}}\,\mathcal{G}^{ik}\,\mathcal{G}^{jl} - a_{ijk}\,a_{lmn}\,\mathcal{G}^{il}\,\mathcal{G}^{jm}\,\mathcal{G}^{kn}$, with $a_{ijk}=\scp{\partial_{ij} g_{\lambda}}{\partial_{k} g_{\lambda}}$.

Therefore, $\bscp{\nabla_{ij}\,g_{\lambda}}{\nabla_{kl}\,g_{\lambda}}$ $\leq$ $ \bscp{\partial_{ij}\,g_{\lambda}}{\partial_{kl}\,g_{\lambda}}\,\mathcal{G}^{ik}\,\mathcal{G}^{jl}$,
from the positive definiteness of $\mathcal{G}^{ij}$ and $\mathcal{G}_{ij}$. Indeed, if we write $W^{ijk\,lmn}=\mathcal{G}^{il}\mathcal{G}^{jm}\mathcal{G}^{kn}$, and if we gather indices $ijk$ and $lmn$ in the two multi-indices\footnote{This can be seen as a relabelling of the $P^3$ combinations of values for $ijk$ into $P^3$ different one-number indices $I$.} $I=(i,j,k)$ and $L=(l,m,n)$, $W^{IL}$ can be seen as a 2-D matrix in $\Rbb^{P^3\times P^3}$. It is then easy to check that the $P^3$ eigenvectors of $W^{IL}$ are given by the $P^3$ combinations of the product of three of the $P$ eigenvectors of $\mathcal{G}^{ij}$, i.e. the covariant vectors $\zeta_i$ respecting the equation $\mathcal{G}^{ij}\,\zeta_j = \mu\,\delta^{ij}\,\zeta_j$ for a certain $\mu=\mu(\zeta)>0$. The matrix $\mathcal{G}^{ij}$ being positive, the eigenvalues of $W^{IL}$ are thus all positive, and $W^{IL}$ is positive. Therefore, $a_IW^{IL}a_{L}\geq 0$ for any tensor $a_I=a_{ijk}$. 
\cqfd

\section{Proof of Lemma \ref{prop:distance-between-maxima-one-orbit-point}}
\label{app:proof-distance-between-maxima-one-orbit-point}

Recall that we use the gradient ascent defined from the optimization function $\phi_r$ such that $\phi_{r+1}(k) = \phi_r(k) + t_r\,\xi_r(k)$, for a sequence of positive step size $t_r$ increasing $S_u$ at each step, and for a step direction
$\xi^i_r(\lambda) \triangleq |\nabla S_u(\phi_r(\lambda))|^{-1} \nabla^i S_u(\phi_r(\lambda))$. From this definition, starting from $k\in\Lambda$, if $\lim_{r\to+\infty}\phi_r(k) = k^\infty\in\Lambda$ exists, then $k^\infty$ is a point where $\nabla^i S_u(k^\infty)=0$ for all $i$, since $S_u(\phi_{r+1}(k)) = S_u(\phi_{r}(k)) + t_r\,|\partial S_u(\phi_r(k))| + O(t_r^2)$.
\medskip

How may the trajectory $\mathcal{T}_k=\{\phi_r(k):r\in\Nbb\}$ contain a point $\lambda'$ satisfying \eqref{eq:max-futur-k-rel} ? Let us write $\gamma_r(s)$ for the geodesic linking $\phi_r(k)$ to $\lambda_M$, and define the \emph{distance function} $\zeta_r=d_{\mathcal{G}}(\lambda_M, \phi_r(k))$. We have thus $\gamma_r(0)=\phi_r(k)$ and $\gamma_r(\zeta_r)=\lambda_M$, where $\lambda_M$ is the global maximum of $S_u$.
\medskip

\noindent\emph{\underline{Case 1.}\quad If $\xi_0^i{\gamma'_0}^j(0)\,\mathcal{G}_{ij}(k) < 0$, i.e. the optimization starts in the wrong direction.} The function $\psi(s)=S_u(\gamma_0(s))$ is twice differentiable over $[0,\zeta_0]$ and for $s$ close to zero, we have $\psi(0)>\psi(s)$ since 
$\psi'(0) = \partial_i S_u(k)\,{\gamma'_0}^i(0)\ = |\nabla S_u(k)|\,\xi_0^i{\gamma'_0}^j(0)\,\mathcal{G}_{ij}(k)\ <\ 0$.

Since $\lambda_M$ is a global maximum of $S_u$, $\psi(0)<\psi(\zeta_0)=S_u(\lambda_M)$. Therefore, there exists a $s^*\in(0,\zeta_0)$ that minimizes $\psi$, i.e. $\psi'(s^*)=0$ with $\psi(s^*)<\psi(0)$. For $\lambda_*=\gamma_0(s^*)$, this implies that $\lambda_*$ is critical since $\psi'(s^*)=\partial_i S_u(\lambda_*){\gamma'_0}^i(s^*)=0$. From Lemma \ref{prop:distance-between-extrema}, 
$ S_u(\lambda_M)-S_u(\lambda_*) \leq \tinv{2}\|u\|^2d_{\mathcal{G}}(\lambda_M,\lambda_*)^2\,(1+\mathcal{K}) < \tinv{2}\|u\|^2d_{\mathcal{G}}(\lambda_M,k)^2\,(1+\mathcal{K})$, since $d_{\mathcal{G}}(\lambda_M,\lambda_*)<d_{\mathcal{G}}(\lambda_0,k)$. Finally, for any $\lambda'\in\mathcal{T}_k$, $S_u(\lambda_M)-S_u(\lambda') \leq S_u(\lambda_M)-S_u(k)$, and $S_u(\lambda_M)-S_u(\lambda') \leq S_u(\lambda_M) - S_u(\lambda_*) \leq  \tinv{2}\|u\|^2\,d_{\mathcal{G}}(\lambda_M,k)^2\,(1+\mathcal{K})$, since $S_u(k) \geq S(\lambda_*)$.
\medskip

\noindent\emph{\underline{Case 2.}\quad If $\xi_0^i{\gamma'_0}^j(0)\,\mathcal{G}_{ij}(k) = 0$.} We have right away ${\gamma'_0}^i(0)\partial_iS_u(k)=0$, and $k$ is a critical point in the direction $\lambda_M$. Lemma \ref{prop:distance-between-extrema} applied on $k$ gives $S_u(\lambda_M)-S_u(k) \leq \tinv{2}\,\|u\|^2\,d_{\mathcal{G}}(\lambda_M,k)^2\,(1+\mathcal{K})$. Since $S_u(\lambda_M)-S_u(\lambda')\leq S_u(\lambda_M)-S_u(k)$ for any $\lambda'\in\mathcal{T}_k$, Equation \eqref{eq:max-futur-k-rel} holds.
\medskip

\noindent\emph{\underline{Case 3.}\quad If $\xi_0^i{\gamma'_0}^j(0)\,\mathcal{G}_{ij}(k) > 0$.} Let us analyze the behavior of the distance function $\zeta_r$.

Let us introduce the function $d_M(\lambda)=d_{\mathcal{G}}(\lambda_M,\lambda)$. As for the Euclidean space, it is easy to prove\footnote{The interested reader will find a proof of this basic differential geometry result in the companion Technical Report \cite{TRLJ0701}.} that $\nabla^id_M(\lambda)=-\gamma^i(0)$ if $\gamma$ is the geodesic linking $\lambda=\gamma(0)$ to $\lambda_M$. Therefore, since $\zeta_{r+1}=d_M(\phi_{r+1}(k))$, a Taylor expansion of $d_M(\lambda)$ around $\lambda=\phi_r(k)$ provides
\begin{equation}
\label{eq:first-order-evolution-zeta_r}
\zeta_{r+1}\ =\ \zeta_r\ -\ t_r\,\xi^i_r(k)\,{\gamma'_r}^j(0)\,\mathcal{G}_{ij}(\phi_r(k))\ +\ O(t_r^2).
\end{equation}
For $r=0$, if $t_0$ is sufficiently small, $\zeta_1<\zeta_0$ and $\zeta_r$ has either a local minima on a particular step $r_m>0$, or it decreases monotically and converges to a value $\zeta_\infty = \lim_{r\to\infty}\zeta_r < \zeta_0$. 

\medskip
\emph{(i) $\zeta$ has a local minima $\zeta_{r_m}<\zeta_0$ on $r_m>0$~:} Then, $\zeta_{r_m+1}>\zeta_{r_m}$ and, using \eqref{eq:first-order-evolution-zeta_r} with some implicit dependences, $\zeta_{r_m+1}-\zeta_{r_m} = -\,t_{r_m}\,\gamma_{r_m}^{\prime\,i}(0)\,\xi^j_{r_m}\,\mathcal{G}_{ij} + O(t_{r_m}^2)$. Therefore, for a sufficiently small step $t_{r_m}$, $\gamma_{r_m}^{\prime\,i}(0)\,\xi^j_{r_m}\,\mathcal{G}_{ij}<0$ and we are in the same hypothesis as \emph{Case 1} with the point $\lambda'=\phi_{r_m}(k)\in\mathcal{T}_k$ instead of $k$. We obtain then
$ S_u(\lambda_M)-S_u(\lambda') \leq \tinv{2}\|u\|^2\,d_{\mathcal{G}}(\lambda_M,\lambda')^2\,(1+\mathcal{K}) < \tinv{2}\|u\|^2\,d_{\mathcal{G}}(\lambda_M,k)^2\,(1+\mathcal{K})$, since $d_{\mathcal{G}}(\lambda_M,\lambda')=\zeta_{r_m}<\zeta_0=d_{\mathcal{G}}(\lambda_M,k)$.
\medskip 

\emph{(ii) If $\zeta_r$ decreases monotically for $r>0$~:} Since $\zeta_r\geq 0$, the limit $\lim_{r\to\infty}\zeta_r$ exists and converges to $\zeta_\infty<\zeta_0$. However, it is not guaranteed that the sequence $\{\phi_r(k)\}$ converges to a point of $\Lambda$. Fortunately, since for all $r>0$, $\phi_r(k)$ remains in the finite volume $V_0=\{\lambda\in\Lambda:d_M(\lambda)\leq d_M(k)\}$, this sequence is bounded in the finite dimensional space $\Lambda$. Therefore, from the Bolzano-Weierstrass theorem on the metric space $(\Lambda,d_{\mathcal{G}}(\cdot,\cdot))$, we can find a convergent subsequence $\{r_i\in\Nbb:r_{i+1}>r_i\}$ such that $\lim_{i\to\infty}\phi_{r_i}(k)=k_\infty\in V_0$. On this point, we will have $\nabla^iS_u(k_\infty)=0$ for all $i$. So, $k_\infty$ is an umbilical point and, from Lemma \ref{prop:distance-between-extrema},
$$
S_u(\lambda_M)-S_u(k_\infty)\ \leq\ \tinv{2}\|u\|^2\,d_{\mathcal{G}}(\lambda_M,k_\infty)^2\,(1+\mathcal{K}).
$$
From now on, we abuse notation and write $\phi_{r_i}(k)=\phi_{i}(k)$. Since $\zeta_\infty^2 = d_{\mathcal{G}}(\lambda_M,k_\infty)^2<d_{\mathcal{G}}(\lambda_M,k)^2=\zeta_0^2$, we can find a $\delta>0$ such that $d_{\mathcal{G}}(\lambda_M, k_\infty)^2\,+\,\delta\,<\,d_{\mathcal{G}}(\lambda_M,k)^2$. Therefore, because $\lim_{i\to\infty}S_u(\phi_{i}(k))=S_u(k^\infty)$ by continuity of $S_u$, and since $S_u(\phi_{i}(k))$ increases monotically with $i$, there exists a $i'>0$ such that $S_u(k^\infty) - S(\phi_{i'}(k))\leq\tinv{2}\|u\|^2\,\delta\,(1+\mathcal{K})$. With $\lambda'=\phi_{i'}(k)\in\mathcal{T}_k$, we finally get $S_u(\lambda_M) - S_u(\lambda') < S_u(\lambda_M)\ -\ S_u(k^\infty)\ +\ \tinv{2}\|u\|^2\,\delta\,(1+\mathcal{K})$, so that $S_u(\lambda_M) - S_u(\lambda') \leq \tinv{2}\|u\|^2\ \big(d_{\mathcal{G}}(\lambda_M,k^\infty)^2+\delta\big)\,(1+\mathcal{K}) < \tinv{2}\|u\|^2\,d_{\mathcal{G}}(\lambda_M,k)^2\,(1+\mathcal{K})$. This gives the result and concludes the proof.\cqfd


\end{document}